\documentclass[12pt]{amsart}
\usepackage{amsmath,amsfonts,amssymb,graphicx,hyperref, enumitem,amsthm,mathabx,amsrefs}
\usepackage{array}

\def\xstrut{\rule[-1.25ex]{0pt}{4ex}}
\newcolumntype{L}{>{\xstrut$}l<{$}}
\newcolumntype{C}{>{\xstrut$}c<{$}}
\setlist[enumerate,1]{label={\arabic*.}}

\newcommand{\WigDName}{\mathcal{D}}
\newcommand{\WigDMat}[1]{\WigDName^{#1}}
\newcommand{\WigDRow}[2]{\WigDMat{#1}_{#2}}

\newcommand{\tildek}[1]{\wtilde{k}\paren{#1}}

\newcommand{\AdSq}{{\mathrm{Ad}}^2}
\newcommand{\SymSq}{{\mathrm{Sym}}^2}
\newcommand{\ExtSq}{{\mathrm{Ext}}^2}

\newcommand{\vpmpm}[1]{v_{_{#1}}}
\newcommand{\chipmpm}[1]{\chi_{_{#1}}}
\newcommand{\bv}{\frak{v}}
\newcommand{\bu}{\frak{u}}

\newcommand{\specmu}{\mathbf{spec}}
\newcommand{\cosmu}{\mathbf{cos}}
\newcommand{\sinmu}{\mathbf{sin}}

\newcommand{\paren}[1]{\ensuremath{\left( #1 \right)}}
\newcommand{\set}[1]{\ensuremath{\left\{ #1 \right\}}}
\newcommand{\abs}[1]{\ensuremath{\left| #1 \right|}}
\newcommand{\norm}[1]{\ensuremath{\left\| #1 \right\|}}
\newcommand{\setdiv}{\,\middle|\,}
\newcommand{\summod}[1]{\ensuremath{\,(\mathrm{mod}\,#1)}}

\newcommand{\piecewise}[1]{\left\{\begin{matrix}#1\end{matrix}\right.}

\newcommand{\If}{\mbox{if }}
\newcommand{\Otherwise}{\mbox{otherwise}}

\renewcommand{\Re}{{\mathop{\mathgroup\symoperators Re}}}
\renewcommand{\Im}{{\mathop{\mathgroup\symoperators Im}}}
\newcommand{\sgn}{{\mathop{\mathgroup\symoperators \,sgn}}}

\newcommand{\Z}{\mathbb{Z}}
\newcommand{\R}{\mathbb{R}}
\newcommand{\N}{\mathbb{N}}
\newcommand{\Q}{\mathbb{Q}}
\newcommand{\C}{\mathbb{C}}

\newcommand{\wbar}[1]{\overline{#1}}
\newcommand{\wtilde}[1]{\widetilde{#1}}

\newcommand{\wcheck}[1]{\widecheck{#1}} 
\newcommand{\e}[1]{e\paren{#1}}
\newcommand{\trans}[1]{{#1}^T}

\DeclareMathOperator*{\res}{res}
\DeclareMathOperator{\diag}{diag}

\theoremstyle{plain} 
\newtheorem{thm}{Theorem}
\newtheorem{cor}[thm]{Corollary}
\newtheorem{lem}[thm]{Lemma}
\newtheorem{prop}[thm]{Proposition}

\title{Kuznetsov, Petersson and Weyl on GL(3), I:\\ The principal series forms.}
\author{Jack Buttcane}
\date{29 January 2017}
\address{Mathematics Department, 244 Mathematics Building, Buffalo, NY 14260, USA}
\email{buttcane@buffalo.edu}
\thanks{During the time of this research, the author was supported by NSF grant DMS-1601919.}

\begin{document}

\begin{abstract}
The Kuznetsov and Petersson trace formulae for $GL(2)$ forms may collectively be derived from Poincar\'e series in the space of Maass forms with weight.
Having already developed the spherical spectral Kuznetsov formula for $GL(3)$, the goal of this series of papers is to derive the spectral Kuznetsov formulae for non-spherical Maass forms and use them to produce the corresponding Weyl laws.
Aside from general interest in new types of automorphic forms, this is a necessary step in the development of a theory of exponential sums on $GL(3)$.
We take the opportunity to demonstrate a sort of minimal method for developing Kuznetsov-type formulae, and produce auxillary results in the form of generalizations of Stade's formula and Kontorovich-Lebedev inversion.
This first paper is limited to the non-spherical principal series forms as there are some significant technical details associated with the generalized principal series forms, which will be handled in a separate paper.
The best analog of this type of form on $GL(2)$ is the forms of weight one which sometimes occur on congruence subgroups.
\end{abstract}

\subjclass[2010]{Primary 11F72; Secondary 11F55}

\maketitle

\section{Introduction}

In the papers \cite{HWI,HWII}, we considered the structure of the $GL(3)$ Maass forms which have non-trivial dependence on $SO(3,\R)$.
We identified three distinct types of cusp forms; in the representation-theoretic language, these are spherical principal series ($d=0$), the non-spherical principal series ($d=1$), and the generalized principle series forms ($d \ge 2$), c.f. \cite[Theorem 3]{HWII} parts 1, 2, and 3, respectively.
To apply analysis on such forms, we must necessarily count them, which leads to the development of Weyl laws.
We consider a type of arithmetically-weighted Weyl law which counts asymptotically the number of forms with Langlands parameters in some region as the volume tends to infinity.
The arithmetic weights are the natural adjoint-square weights arising in a type of Kuznetsov trace formula attached to each type of cusp form; for unweighted Weyl laws, see \cite[Theorem 1.2]{Miller01} and \cite[Theorem 0.1]{Muller01}, and for the weighted case, up to a constant, see \cite[Theorem 1]{Val01}.

On $GL(2)$, one may collectively treat the Kuznetsov and Petersson trace formulae as expressions for Fourier coefficients of Poincar\'e series at the minimal $K$-types of Maass forms with weight.
In this way, the Petersson trace formula can be thought of as the Kuznetsov trace formula at a point in the spectrum.\footnote{The author would like to thank Prof. William Duke for this perspective.}
This leads to three possible spectral Kuznetsov formulae and corresponding Weyl laws on $GL(3)$, one for each of the three minimal $K$-types described above.
The Kuznetsov formula for spherical principal series forms was previously considered in the paper \cite{SpectralKuz}, the first paper of this series will consider the non-spherical principal series, which is Theorem \ref{thm:Kuznetsov}, and the second paper will treat the generalized principal series.

The Weyl law for these non-spherical principal series forms is as follows:
Let $\mathcal{S}^1_3$ be an orthonormal basis of vector-valued cusp forms attached to the 3-dimensional representation of $SO(3,\R)$ with $V$-characters $\chipmpm{+-}$ (see \cite[section 2.2.2]{HWI}).
Denote the spectral parameters of such a Maass form $\varphi$ by
\[ \mu_\varphi \in \mathfrak{a}_\C^* := \set{\mu\in\C^3 \setdiv \mu_1+\mu_2+\mu_3=0}. \]
The condition on the $V$-character simply restricts the allowed permutations of $\mu_\varphi$ (see \cite[Theorem 3.2]{HWII}); as a Langlands quotient, these have the form
\begin{align}
\label{eq:LangQuot}
	J(GL(3,\R), P_{1,1,1}, \abs{\cdot}^{\mu_1}\sgn, \abs{\cdot}^{\mu_2}\sgn, \abs{\cdot}^{\mu_3}).
\end{align}
\begin{thm}
\label{thm:HeckeWeyl}
	Take a bounded set
	\[ \Omega \subset i\mathfrak{a}_\R^* := \set{\mu\in\mathfrak{a}_\C^*\setdiv\Re(\mu)=0}, \]
	which is symmetric under $\mu_1 \leftrightarrow \mu_2$ and whose boundary has Minkowski dimension at most 1.
	Also let $\mu'\in i\mathfrak{a}_\R^*$ and $T > M > 1$ be large parameters, then if $\mathcal{S}^1_3$ consists of Hecke eigenforms,
	\begin{align*}
		\sum_{\mu_\varphi \in T\Omega} \frac{1}{L(1,\AdSq \varphi)} =& \frac{3}{2\pi} \int_{T\Omega} \specmu^1(\mu) d\mu+O\paren{T^{4+\epsilon}}, \\
		\sum_{\norm{\mu_\varphi-T \mu'} < M} \frac{1}{L(1,\AdSq \varphi)} =& \frac{3}{2\pi} \int_{\norm{\mu-T \mu'} < M} \specmu^1(\mu) d\mu+O\paren{T^{3+\epsilon}M^{1+\epsilon}}.
	\end{align*}
\end{thm}
The spectral measure $\specmu^1(\mu) d\mu$ here is given by $d\mu=d\mu_1 \, d\mu_2$ and
\begin{align*}
	\specmu^1(\mu) =& \frac{1}{64\pi^4} (\mu_1-\mu_2) (\mu_1-\mu_3)(\mu_2-\mu_3) \cot\frac{\pi}{2}(\mu_1-\mu_3) \cot\frac{\pi}{2}(\mu_2-\mu_3) \tan\frac{\pi}{2}(\mu_1-\mu_2).
\end{align*}
The spectral integrals are of size $\asymp T^5$ and $\asymp T^3 M^2$ provided $\Omega$ contains an open ball and $\specmu^1(\mu')\ne 0$; we expect that the error terms $T^4$ and $T^3 M$ are best-possible due to the sharp cut-off.
We use the Minkowski dimension \cite{Falconer} assumption here as it is exactly the correct hypothesis for the proof, but essentially we are just asking that the boundary not be a fractal.

A word on the history of non-spherical forms for $SL(3,\Z)$:
It has been known for some time \cite{GelbJacq} that symmetric squares of $SL(2,\Z)$ cusp forms are cusp forms for $SL(3,\Z)$, and it is clear that the symmetric square of a holomorphic modular form cannot be spherical; in \cite[section 6.4]{MeMiller}, we note that the symmetric square of a $GL(2)$ form of weight $k$ has weight $d=2k-1$.
Stephen D. Miller has communicated to the author an argument of the late Jonathan Rogawski to the effect that the only $GL(2)$ forms whose symmetric squares have full level are themselves full-level, but possibly twisted by quadratic global Hecke characters.
These two arguments imply the symmetric squares miss the $d=1$ forms, because there are no weight-one forms of full level for $SL(2,\Z)$, by parity considerations.

The development of cohomological forms for $GL(3)$ gives another possible cache to investigate, but Ash and Pollack \cite{Ash01} have shown for small weight that the only cohomological forms of full level are already given by the symmetric-square construction, and they conjecture that this holds for all weight.
M\"uller, however, has given a very general Weyl law \cite[Theorem 0.1]{Muller01}, which confirms the existence of many weight-one forms for $SL(3,\Z)$.
Combined with the knowledge that the spherical forms have no image in the weight-one forms, this confirms the existence of non-spherical forms which do not arise from the symmetric-square construction.

Theorem \ref{thm:HeckeWeyl} gives an arithmetically-weighted version of M\"uller's Weyl law, with a strong error term, again confirming the existence of many non-spherical forms which are not symmetric squares.
In the next paper in this series \cite{WeylII}, we will demonstrate the existence of many cusp forms of minimal weight $d \ge 2$, and these cannot be counted by M\"uller's Weyl law, as they are lost among the lifts of the spherical and weight-one forms.

We also use this opportunity to showcase what appears to be the minimal method for obtaining such Kuznetsov formulae.
The main technical difficulties are the appropriate generalization -- Theorem \ref{thm:StadesFormula} -- of Stade's formula \cite[Theorem 1.1]{Stade02} and the analytic continuation of a certain unpleasant integral (section \ref{sect:LongEleXCont}).
Kontorovich-Lebedev inversion -- Theorem \ref{thm:KontLebedev} -- can be shown directly from Stade's formula as in \cite{GoldKont}, and having already solved the requisite differential equations, the analytic continuation step completes the Kuznetsov formula using the arguments of \cite{SpectralKuz}.

Stade's formula is an evaluation of the Archimedian local zeta integral for Rankin-Selberg on $GL(3)\times GL(3)$ (here in the ramified case), which is of independent interest, see \cite{HiIshMiya}.

\section{Results}
In this paper, we deal primarily with a certain completion of the Jacquet-Whittaker function.
For $\mu\in\mathfrak{a}_\C^*$, $s\in\C^2$, $\beta,\eta\in\Z^3$, define
\begin{align}
\label{eq:GDef}
	\wtilde{G}(1,\beta,\eta, s,\mu) =& \frac{\prod_{i=1}^3 \Gamma\paren{\frac{\beta_i+s_1-\mu_i}{2}} \Gamma\paren{\frac{\eta_i+s_2+\mu_i}{2}}}{\Gamma\paren{\frac{s_1+s_2+\sum_i (\beta_i+\eta_i)-2}{2}}},
\end{align}
and for $\abs{m'}\le 1$, writing $m'=\varepsilon m$ with $\varepsilon = \pm 1$ and $0 \le m \le 1$, set
\begin{align*}
	G^1_{m'}(s,\mu) =& \sqrt{\binom{2}{1+m}} \sum_{\ell=0}^{m} \varepsilon^\ell \binom{m}{\ell} \wtilde{G}(1,(1-m,1-m,m),(\ell,\ell,1-\ell),s,\mu),
\end{align*}
and take $G^1(s,\mu)$ to be the row vector with coordinates $G^1_{m'}(s,\mu)$, indexing from the central entry: $G^1 = (G^1_{-1},G^1_0,G^1_1)$.
Then we define the completed Whittaker function for the non-spherical principal series forms as
\begin{align*}
	W^{1*}(y, \mu) =& \frac{1}{4\pi^2} \int_{\Re(s)=\mathfrak{s}} (\pi y_1)^{1-s_1} (\pi y_2)^{1-s_2} G^1\paren{s,\mu} \frac{ds}{(2\pi i)^2},
\end{align*}
c.f. \cite[Theorem 5]{HWII}, where $\mathfrak{s}$ may be taken anywhere to the right of the poles of the integrand.
We extend to $G=PSL(3,\R)$ by the Iwasawa decomposition
\[ W^{1*}(xyk, \mu) = \psi_{1,1}(x) W^{1*}(y, \mu) \WigDMat{1}(k), \]
where $\WigDMat{1}:K\to GL(3,\C)$ is the 3-dimensional Wigner $\WigDName$-matrix (see \cite[section 2.2]{HWI}).

Our first step is to generalize Stade's formula to these Whittaker functions, which we will use twice, in different contexts.
\begin{thm}
\label{thm:StadesFormula}
	Define
	\[ \Psi^1=\Psi^1(\mu,\mu',t) = \int_{Y^+} W^{1*}(y,\mu) \trans{W^{1*}(y,\mu')} (y_1^2 y_2)^t dy, \]
	with $(Y^+, dy)$ as in \cite[section 2.1]{HWI}, and set
	\[ c^1_{i,j} = \piecewise{1 & \If i=3\ne j \text{ or } j=3\ne i, \\ 0 & \Otherwise,} \]
	then
	\begin{align*}
		\Psi^1 =& \frac{1}{2\pi^{3t} \Gamma\paren{\tfrac{3t}{2}}} \prod_{i,j} \Gamma\paren{\tfrac{c^1_{i,j}+t+\mu_i'+\mu_j}{2}}.
	\end{align*}
\end{thm}
One could more generally consider inserting some $\WigDMat{1}(v), v\in V$ between the Whittaker functions, but the current formula is sufficient for our purposes.

Theorem \ref{thm:StadesFormula} confirms a generalization of the conjecture of Bump \cite[section 2.6]{Bump02} that $\Psi^1$ gives the gamma factors of the Rankin-Selberg convolution $L$-function; it has direct, immediate implications for such $L$-functions:
\begin{cor}
\label{cor:RSLfunc}
	Let $\varphi,\varphi'\in \mathcal{S}^1_3$ be Hecke eigenforms.
	Then the completed Rankin-Selberg $L$-function
	\[ \Lambda(s,\varphi \times \wbar{\varphi'}) = 2\Gamma_\R(3s) \Psi^1(\mu_{\varphi}, \wbar{\mu_{\varphi'}},s) L(s,\varphi \times \wbar{\varphi'}), \qquad \Gamma_\R(s) = \pi^{-s/2} \Gamma(\tfrac{s}{2}) \]
	has functional equation 
	\[ \Lambda(s,\varphi \times \wbar{\varphi'}) = \Lambda(1-s,\wbar{\varphi} \times \varphi'), \]
	and is entire except for possible simple poles at $s=\frac{1}{2} \pm \frac{1}{2}$ with residues
	\[ \pm \frac{2}{3} \int_{\Gamma\backslash G} \varphi(g) \wbar{\trans{\varphi'(g)}} dg. \]
\end{cor}
There is a considerable history for results of this type, and some further discussion can be found in \cite{Fried01} or \cite{Stade02}.

We will require a non-trivial bound on the spectral parameters; this was essentially proved in \cite{KS01}:
\begin{thm}[Kim, Sarnak]
\label{thm:KS}
	The spectral parameters $\mu$ of a cusp form of weight one satisfy $\abs{\Re(\mu_j)} < \frac{5}{14}$.
\end{thm}
Even though \cite[Proposition 1 of appendix 2]{KS01} (and \cite[Theorem 1.2]{LRS} from which the method derives) assumes the cusp form in question is spherical, it really only relies on the existence of the gamma factors 
\[ \Gamma_\R\paren{s+2\mu_1} \Gamma_\R\paren{s+2\mu_2} \Gamma_\R\paren{s+2\mu_3} \]
in the completed symmetric-square $L$-function, and this continues to hold for weight-one forms, as well.
We will discuss this and Corollary \ref{cor:RSLfunc} further in section \ref{sect:RS}.

The first application of Stade's formula in the derivation of the weight-one Kuznetsov formula is to generalize Kontorovich-Lebedev inversion for our choice of Whittaker functions, using the method of Goldfeld and Kontorovich \cite{GoldKont}.
\begin{thm}
\label{thm:KontLebedev}
	Set
	\begin{align*}
		\frac{1}{\sinmu^1(\mu)} =& \frac{(2\pi i)^2}{3} \lim_{t\to0} t^2 \Psi^1(\mu,-\mu,t).
	\end{align*}
	For $f:Y^+\to\C^3$ define
	\[ f^\sharp(\mu) = \int_{Y^+} f(y) \wbar{\trans{W^{1*}(y,\mu)}} dy, \]
	and for $F:\mathfrak{a}_\C^*\to\C$, define
	\[ F^\flat(y) = \int_{\Re(\mu)=0} F(\mu) W^{1*}(y,\mu) \sinmu^1(\mu)\,d\mu. \]
	If $F(\mu)$ is holomorphic and Schwartz-class on a tube domain $\set{\mu\setdiv\abs{\Re(\mu_i)}<\delta}$ for some $\delta > 0$ and invariant under $\mu_1 \leftrightarrow \mu_2$, then
	\[ (F^\flat)^\sharp(\mu) = F(\mu). \]
\end{thm}
Here we are making no claim as to the image of $F \mapsto F^\flat$ beyond the necessary convergence.
Note that one may remove the invariance condition for $F$, and the conclusion becomes
\[ 2(F^\flat)^\sharp(\mu) = F(\mu)+F(\mu_2,\mu_1,\mu_3). \]
The Kontorovich-Lebedev spectral measure may also be written as
\[ \sinmu^1(\mu) = \frac{(\mu_1-\mu_2) \cos\frac{\pi}{2}(\mu_1-\mu_3) \cos\frac{\pi}{2}(\mu_2-\mu_3) \sin\frac{\pi}{2}(\mu_1-\mu_2)}{16\pi^5}. \]

We then consider a Fourier coefficient of a Poincar\'e series defined by summing an inverse Whittaker transform.
\begin{align}
\label{eq:MainPoincare}
	P_m(g,F) =& \sum_{\gamma\in U(\Z)\backslash\Gamma} \int_{\Re(\mu)=0} F(\mu) W^{1*}(\wtilde{m} \gamma g,\mu) \sinmu^1(\mu)\,d\mu, \qquad \wtilde{m}=\diag(m_1 m_2,m_1,1).
\end{align}
The usual spectral expansion and Bruhat decomposition give our pre-Kuznetsov trace formula, and the second application of Stade's formula completes the spectral side and the trivial term of the arithmetic/geometric side.
The spectral measure occuring in the trivial term can be thought of as Stade's formula at 1 divided by Stade's formula at 0.

The paper \cite{Spectral Kuz} was made slightly more delicate than necessary by the inclusion of the final $y$-integral, essentially leap-frogging the pre-Kuznetsov formula.
Here, we derive the Kuznetsov kernel functions in two steps; first applying uniqueness of the Kuznetsov kernel functions, as in section \ref{sect:KernelDiffEqs}, and then applying Stade's formula to complete the analysis.
Having already solved the requisite differential equations, the key ingredient is the analytic continuation of a certain unpleasant integral, but this is roughly similar to the spherical case, so not much additional work is needed.

For the Kuznetsov kernel functions, using the power series solutions $J_{w_l}(y,\mu)$ as in section \ref{sect:KernelDiffEqs}, we set
\begin{align*}
	J^1_{w_l}(y,\mu) :=& \varepsilon_2 J_{w_l}(y,\mu)+\varepsilon_1 J_{w_l}(y,\mu^{w_4})+\varepsilon_1 \varepsilon_2 J_{w_l}(y,\mu^{w_5}), \\
	J^1_{w_4}(y,\mu) :=& -\sin\frac{\pi}{2}(\mu_1-\mu_2) J_{w_4}(y,\mu) - i\varepsilon_1 \cos\frac{\pi}{2}(\mu_1-\mu_3) J_{w_4}(y,\mu^{w_4}) \\
	& \qquad + i\varepsilon_1 \cos\frac{\pi}{2}(\mu_2-\mu_3) J_{w_4}(y,\mu^{w_5}),
\end{align*}
where $\varepsilon=\sgn(y)$ (note the arguments of the $J_w(y,\mu)$ functions are still the \textbf{signed} $y$), and define
\begin{align}
\label{eq:K1I}
	K^1_I(y,\mu) :=& 1, \\
\label{eq:K1w4}
	K^1_{w_4}(y,\mu) :=& \frac{1}{8\pi} \frac{J^1_{w_4}(y,\mu)}{\cos\frac{\pi}{2}(\mu_1-\mu_3)\cos\frac{\pi}{2}(\mu_2-\mu_3)\sin\frac{\pi}{2}(\mu_1-\mu_2)}, \\
\label{eq:K1w5}
	K^1_{w_5}(y,\mu) :=& K^1_{w_4}((-y_2,y_1),-\mu), \\
\label{eq:K1wl}
	K^1_{w_l}(y,\mu) :=& -\frac{1}{16\pi} \frac{J^1_{w_l}(y,\mu)-J^1_{w_l}(y,\mu^{w_2})}{\cos\frac{\pi}{2}(\mu_1-\mu_3)\cos\frac{\pi}{2}(\mu_2-\mu_3)\sin\frac{\pi}{2}(\mu_1-\mu_2)}.
\end{align}

Then we define the integral transforms
\begin{align*}
	H_w(F; y) =& \frac{1}{\abs{y_1 y_2}} \int_{\Re(\mu)=0} F(\mu) K^1_w(y, \mu) \specmu^1(\mu) d\mu,
\end{align*}
with
\begin{align*}
	\frac{1}{\cosmu^1(\mu)} :=& \Psi^1(\mu,-\mu,1) = \frac{\pi (\mu_1-\mu_3)(\mu_2-\mu_3)}{4 \sin\frac{\pi}{2}(\mu_1-\mu_3) \sin\frac{\pi}{2}(\mu_2-\mu_3) \cos\frac{\pi}{2}(\mu_1-\mu_2)}, \\
	\specmu^1(\mu) :=& \frac{\sinmu^1(\mu)}{\cosmu^1(\mu)}.
\end{align*}

For a cusp form $\varphi \in \mathcal{S}^1_3$, we denote its Fourier-Whittaker coefficients as $\rho_\varphi^*(m)$, $m\in\Z^2$ and its Langlands parameters as $\mu_\varphi$, see \cite{HWII}.
Similarly, let $\mathcal{S}^1_2$ be a basis of $SL(2,\Z)$ spherical Maass cusp forms.
For such a cusp form $\phi$ we denote the Fourier-Whittaker coefficients of the maximal parabolic Eisenstein series attached to $\phi$ with the additional spectral parameter $\mu_1$ as $\rho_\phi^*(m;\mu_1)$, see \cite[section 5]{HWI}.

The $SL(3,\Z)$ Kloosterman sum attached to the characters $\psi_m$ and $\psi_n$ and the Weyl element $w$ with moduli $c\in\N^2$ is denoted $S_w(\psi_m,\psi_n,c)$, see section \ref{sect:Kloos}.

\begin{thm}
\label{thm:Kuznetsov}
	Let $F(\mu)$ be Schwartz-class, holomorphic on $\set{\mu|\abs{\Re(\mu_i)} < \frac{1}{2}+\delta}$ for some $\delta > 0$, and invariant under $\mu_1 \leftrightarrow \mu_2$.
	Suppose $F(\mu)=0$ whenever $\mu_1-\mu_2=\pm1$, then for $m,n\in\Z^2$ with $m_1 m_2 n_1 n_2 \ne 0$,
	\[ \mathcal{C}+\mathcal{E}=\mathcal{K}_I+\mathcal{K}_4+\mathcal{K}_5+\mathcal{K}_l, \]
	where
	\begin{align*}
		\mathcal{C} =& \sum_{\varphi \in \mathcal{S}_3^1} F(\mu_\varphi) \frac{\wbar{\rho_\varphi^*(m)} \rho_\varphi^*(n)}{\cosmu^1(\mu_\varphi)}, \\
		\mathcal{E} =&\frac{1}{2\pi i} \sum_{\substack{\phi \in \mathcal{S}_2^1\\ \phi\text{ odd}}} \int_{\Re(\mu_1)=0} F(\mu_1+\mu_\phi,\mu_1-\mu_\phi,-2\mu_1) \frac{\wbar{\rho_\phi^*(m;\mu_1)} \rho_\phi^*(n;\mu_1)}{\cosmu^1(\mu_1+\mu_\phi,\mu_1-\mu_\phi,-2\mu_1)} d\mu_1,
	\end{align*}
	\begin{align*}
		\mathcal{K}_I =& \delta_{\substack{\abs{m_1}=\abs{n_1}\\ \abs{m_2}=\abs{n_2}}} H_I(F;(1, 1)) \nonumber \\
		\mathcal{K}_4 =& \sum_{\varepsilon\in\set{\pm1}^2} \sum_{\substack{c_1,c_2\in\N\\ \varepsilon_1 m_2 c_1=n_1 c_2^2}} \frac{S_{w_4}(\psi_m,\psi_{\varepsilon n},c)}{c_1 c_2} H_{w_4}\paren{F; \paren{\varepsilon_1 \varepsilon_2 \tfrac{m_1 m_2^2n_2}{c_2^3n_1},1}}, \\
		\mathcal{K}_5 =& \sum_{\varepsilon\in\set{\pm1}^2} \sum_{\substack{c_1,c_2\in\N\\ \varepsilon_2 m_1 c_2= n_2 c_1^2}} \frac{S_{w_5}(\psi_m,\psi_{\varepsilon n},c)}{c_1 c_2} H_{w_5}\paren{F; \paren{1,\varepsilon_1 \varepsilon_2 \tfrac{m_1^2 m_2 n_1}{c_1^3 n_2}}}, \\
		\mathcal{K}_l =& \sum_{\varepsilon\in\set{\pm1}^2} \sum_{c_1,c_2\in\N} \frac{S_{w_l}(\psi_m,\psi_{\varepsilon n},c)}{c_1 c_2} H_{w_l}\paren{F; \paren{\varepsilon_2 \tfrac{m_1 n_2 c_2}{c_1^2}, \varepsilon_1 \tfrac{m_2 n_1 c_1}{c_2^2}}}.
	\end{align*}
\end{thm}

It is interesting to note that the existence of the Kontorovich-Lebedev inversion, the need for the index integral and the existence of a pole in Stade's formula at 0 roughly coincide.
Otherwise, we could consider absolutely convergent Poincar\'e series formed from Whittaker functions at individual points in the spectrum, similar to the Petersson trace formulae on $GL(2)$.
This will be reflected in the formulae for the generalized principal series representations to be considered in the second part, where the index integral will be one-dimensional.

We may compare the linear combinations of the power series solutions to the spherical case:
\begin{prop}
For $K_{w_l}^{\pm\pm}$, $K_{w_l}$ and $\sinmu(\mu)$ as in section \ref{sect:OldMBIntegrals},
\begin{enumerate}

\item For $y_1 < 0 < y_2$, \(\displaystyle J^1_{w_l}(y,\mu)-J^1_{w_l}(y,\mu^{w_2}) = K_{w_l}^{-+}(y,\mu)+K_{w_l}^{-+}(y,\mu^{w_3})+K_{w_l}^{-+}(y,\mu^{w_l})\),

\item For $y_2 < 0 < y_1$, \(\displaystyle J^1_{w_l}(y,\mu)-J^1_{w_l}(y,\mu^{w_2}) = -K_{w_l}^{+-}(y,\mu)+K_{w_l}^{+-}(y,\mu^{w_2})-K_{w_l}^{+-}(y,\mu^{w_l})\),

\item For $y_1,y_2 < 0$, \(\displaystyle J^1_{w_l}(y,\mu)-J^1_{w_l}(y,\mu^{w_2}) = -K_{w_l}^{--}(y,\mu)+K_{w_l}^{--}(y,\mu^{w_2})-K_{w_l}^{--}(y,\mu^{w_3})\),

\item For $y_1,y_2 > 0$, \(\displaystyle J^1_{w_l}(y,\mu)-J^1_{w_l}(y,\mu^{w_2}) = -\tfrac{32}{\pi^3} \sinmu(\mu) K_{w_l}(y,\mu)\).

\end{enumerate}
\end{prop}
We will use this proposition for the Mellin-Barnes integrals given in \cite[Theorem 2]{SpectralKuz}.
Note that the $K_{w_l}^{\pm\pm}$ functions defined in \cite{Subconv} also include the factor $1/\sinmu(\mu)$, while the functions in \cite{SpectralKuz} do not.
There is a small point that there is no agreement in the signs of the above linear combinations with the spherical case; for the latter, the sign on $K_{w_l}^{\pm\pm}(y,\mu^w)$ is necessarily the number of transpositions $w_2$, $w_3$ in the permutation $w$.

The $K_{w_4}^1$ function is new, but it is easy to verify that
\begin{prop}
	For $y_1\ne 0$,
	\begin{align*}
		K_{w_4}^1((y_1,1),\mu) =& \frac{1}{8\pi^3} \int_{-i\infty}^{i\infty}\abs{8\pi^3y_1}^{1-s}\Gamma(s-\mu_1)\Gamma(s-\mu_2)\Gamma(s-\mu_3) \biggl(\exp\paren{i\frac{3\pi}{2}\varepsilon s}\\
		&-\exp\paren{-i\frac{\pi}{2}\varepsilon(s+2\mu_1)}-\exp\paren{-i\frac{\pi}{2}\varepsilon(s+2\mu_2)}+\exp\paren{-i\frac{\pi}{2}\varepsilon(s+2\mu_3)}\biggr) \frac{ds}{2\pi i},
	\end{align*}
	where $\varepsilon=\sgn(y_1)$.
\end{prop}

We state a technical version of the Weyl law with both analytic and arithmetic weights:
\begin{thm}
\label{thm:TechnicalWeyl}
	Let $F(\mu)$ and $\delta$ be as in Theorem \ref{thm:Kuznetsov}, and define
	\[ E_F(s,t) = \int_{\Re(\mu)=s} \paren{\sum_{w\in W} \abs{F(\mu^w)}} \paren{\sum_{w\in W} (1+\abs{\mu^w_1-\mu^w_2})^{t_1+\epsilon} (1+\abs{\mu^w_2-\mu^w_3})^{t_2+\epsilon}} \abs{d\mu}, \]
	\[ E_F^* = \sum_{\substack{\phi \in \mathcal{S}_2^1\\ \phi\text{ odd}}} \int_{\Re(\mu_1)=0} \abs{F(\mu_1+\mu_\phi,\mu_1-\mu_\phi,-2\mu_1)} \paren{1+\abs{\mu_1}}^\epsilon\paren{1+\abs{\mu_\phi}}^\epsilon \abs{d\mu_1}. \]	Then we have
	\begin{align*}
		\sum_{\varphi \in \mathcal{S}_3^1} F(\mu_\varphi) \frac{\abs{\rho_\varphi^*(1)}^2}{\cosmu^1(\mu_\varphi)} = \int_{\Re(\mu)=0} F(\mu) \specmu^1(\mu) d\mu+O\paren{E_F^*+\sum_{i=1}^4 E_F(s_i, t_i)},
	\end{align*}
	with the parameters
	\begin{center}
	\begin{tabular}{L|C|C|C|C}
	i & 1 & 2 & 3 & 4 \\
	\hline
	s_i & (0,0,0) & (\frac{1}{4}+2\eta,\frac{1}{4}+2\eta,-\frac{1}{2}-4\eta) & (\frac{1}{2}+4\eta,-\frac{1}{4}-2\eta,-\frac{1}{4}-2\eta) & (\frac{1}{2}+\eta,0,-\frac{1}{2}-\eta) \\
	\hline
	t_i & (0,\frac{1}{2}) & (-\frac{1}{4},\frac{3}{4}) & (-\frac{1}{4},\frac{3}{4}) & (-\frac{1}{2},0),
	\end{tabular}
	\end{center}
	where we assume $0 < 4\eta < \delta$ satisfies $\eta = O(\epsilon)$.
\end{thm}
We have made no attempt to optimize the error terms.

By a careful choice of test function, we will remove the analytic weights:
\begin{cor}
\label{cor:Weyl}
	For $\Omega$,  $\mu'$ and $T$ and $M$ as in Theorem \ref{thm:HeckeWeyl},
	\begin{enumerate}
	\item \(\displaystyle \sum_{\mu_\varphi \in T\Omega} \frac{\abs{\rho_\varphi^*(1)}^2}{\cosmu^1(\mu_\varphi)} = \int_{T\Omega} \specmu^1(\mu) d\mu+O\paren{T^{4+\epsilon}}\),
	\item \(\displaystyle \sum_{\norm{\mu_\varphi-T \mu'} < M} \frac{\abs{\rho_\varphi^*(1)}^2}{\cosmu^1(\mu_\varphi)} = \int_{\norm{\mu-T \mu'} < M} \specmu^1(\mu) d\mu+O\paren{T^{3+\epsilon}M^{1+\epsilon}}\).
	\end{enumerate}
\end{cor}

We finish by giving the Kuznetsov formula on Hecke eigenvalues, which follows by a Rankin-Selberg argument using Stade's formula for a third time.
\begin{thm}
\label{thm:HeckeKuznetsov}
	When the bases of Theorem \ref{thm:Kuznetsov} are taken to be Hecke eigenfunctions, the left-hand side may be written as
	\begin{align*}
		\mathcal{C} =& \frac{2\pi}{3} \sum_{\varphi \in \mathcal{S}_3^1} F(\mu_\varphi) \frac{\wbar{\lambda_\varphi(m)} \lambda_\varphi(n)}{L(1,\AdSq \varphi)}, \\
		\mathcal{E} =& \frac{2\pi}{2\pi i} \sum_{\substack{\phi \in \mathcal{S}_2^1\\ \phi\text{ odd}}} \int_{\Re(\mu_1)=0} \frac{F(\mu_1+\mu_\phi,\mu_1-\mu_\phi,-2\mu_1) \wbar{\lambda_\phi(m,\mu_1)} \lambda_\phi(n,\mu_1)}{L(\phi,1+3\mu_1) L(\phi,1-3\mu_1) L(1,\AdSq \phi)} d\mu_1, \nonumber
	\end{align*}
	where $\lambda_\varphi(m)$ and $\lambda_\phi(m,\mu_1)$ as in \cite[(5.13)]{HWI} are the the Hecke eigenvalues of the associated forms.
\end{thm}
For $\varphi \in \mathcal{S}_3^1$, $L(1,\AdSq \varphi)$ is the residue at $s=1$ of the Rankin-Selberg $L$-function $L(s,\varphi\times\wbar{\varphi})$, see \eqref{eq:AdSqFactor}.

The main theorem now follows.

\section{Acknowledgements}
The author would like to thank Lior Silberman for pointing out the full generality of M\"uller's Weyl law, and Valentin Blomer for pointing out the applicability of \cite{KS01} and \cite{LRS}.
The author would also like to thank the referee for their helpful comments.

\section{Background}
\label{sect:Background}
We generally follow the notation of \cite{HWI,HWII} and \cite{SpectralKuz}; a quick recollection of much of the notation can be found in \cite[section 2]{HWII}.
There are two main discrepancies between the two notations:
First, the description of the Weyl group given in \cite{HWI} was chosen to detach the elements of the $V$ group, so we use that description over the one in \cite{SpectralKuz}.
Second, the description of the Laplacian in \cite{HWII} has the sign changed to make it a positive operator; we use that description over the one in \cite{SpectralKuz}.

There are two effects of changing the long Weyl element $w_l$:
For the Kloosterman sums $S_{w_l}(\psi_m,\psi_n,c)$, this effectively replaces $\psi_n \mapsto \psi_n^{\vpmpm{--}} = \psi_{-n}$.
The second change has to do with solving the differential equations, but is directly offset by reintroducing the conjugation $\wbar{\psi_{1,1}(x)}$ in \eqref{eq:HwTildeDef}; conversely, the conjugation will affect the $w_4$ and $w_5$ differential equations, but not their Kloosterman sums.
We repeat the power-series solutions in the current notation in section \ref{sect:KernelDiffEqs}, below.

Since we track eigenvalues of the Laplacian in terms of the eigenvalue of the power function, the leading sign of the Laplacian is essentially irrelevant.

Let $\bv^d_j=(\bv^d_{j,-d},\ldots,\bv^d_{j,d})$ be the $(2d+1)$-dimensional row vector with entries
\begin{align}
\label{eq:vdef}
	\bv^d_{j,m'} = \delta_{m'=j},
\end{align}
and set
\begin{align}
\label{eq:udef}
	\bu^{d,\pm}_j = \frac{1}{2}(\bv^d_j \pm (-1)^d \bv^d_{-j}).
\end{align}
These vectors satisfy the parity conditions
\begin{align}
\label{eq:uSigns}
	\bu^{1,-}_0 \WigDMat{1}(\vpmpm{\varepsilon_1,\varepsilon_2})=&\varepsilon_2 \bu^{1,-}_0, &
	\bu^{1,-}_1 \WigDMat{1}(\vpmpm{\varepsilon_1,\varepsilon_2})=&\varepsilon_1 \varepsilon_2 \bu^{1,-}_1, &
	\bu^{1,+}_1 \WigDMat{1}(\vpmpm{\varepsilon_1,\varepsilon_2})=&\varepsilon_1 \bu^{1,+}_1,
\end{align}
where $\vpmpm{\varepsilon_1,\varepsilon_2}=\diag\set{\varepsilon_1,\varepsilon_1\varepsilon_2,\varepsilon_2}$ is an element of the group $V$ of diagonal orthogonal matrices as in \cite[section 2.1]{HWI}.

\subsection{The Whittaker functions}
\label{sect:WhittFuncs}
Define
\begin{align*}
	c_W(I)=&\sqrt{2}&\alpha(I)=&(0,1,1)&\bu(I)=\bu^{1,-}_0\\
	c_W(w_4)=&2&\alpha(w_4)=&(1,1,0)&\bu(w_4)=\bu^{1,+}_1\\
	c_W(w_5)=&-2&\alpha(w_5)=&(1,0,1)&\bu(w_5)=\bu^{1,-}_1,
\end{align*}
and extend by left $w_2$-invariance, i.e. $c_W(w_2 w)=c_W(w)$, etc.
Using 
\[ \Lambda_\alpha(\mu) = \pi^{-\frac{3}{2}+\mu_3-\mu_1} \Gamma\paren{\tfrac{1+\alpha_1+\mu_1-\mu_2}{2}} \Gamma\paren{\tfrac{1+\alpha_2+\mu_1-\mu_3}{2}} \Gamma\paren{\tfrac{1+\alpha_3+\mu_2-\mu_3}{2}} \]
as in \cite[Theorem 5]{HWII}, the Mellin-Barnes integral representations of the Jacquet-Whittaker functions attached to $\WigDMat{1}$ read
\begin{align}
\label{eq:Whitt1FEs}
	W^{1*}(g, \mu) =& c_W(w) \Lambda_{\alpha(w)}(\mu^w) \bu(w) W^1(g, \mu^w, \psi_{1,1}).
\end{align}

From these functional equations, we may deduce (see \cite[section 3.5]{HWI})
\begin{align}
\label{eq:Whitty1Asymp}
	W^{1*}(y,\mu) \sim& \sum_{w\in W_3} c_W(w) p_{\rho+\mu^{w w_l}}(y) \Lambda_{\alpha(w)}(\mu^w) \bu(w) W^1(I,\mu^w,\psi_{0,y_2}), \qquad y_1 \to 0,
\end{align}
when the components of $\mu$ are distinct and on a neighborhood of $\Re(\mu)=0$, say.

Because the argument that follows depends so crucially on the values of these constants, we describe one means of verification:
The first-term asymptotics of the Whittaker functions can be computed directly from the Mellin-Barnes integral, i.e. \eqref{eq:WhittResExpand}, or via the general computation of \cite[section 3.5]{HWI}, or even in two steps using a combination of both via \eqref{eq:Whitty1Asymp} and the asymptotics of the classical Whittaker function:
\begin{align}
\label{eq:ClassWhittAsymp}
	W_{\alpha,\beta}(y) \sim y^{\frac{1}{2}+\beta} \frac{\Gamma\paren{-2\beta}}{\Gamma\paren{\frac{1}{2}-\alpha-\beta}}+y^{\frac{1}{2}-\beta} \frac{\Gamma\paren{2\beta}}{\Gamma\paren{\frac{1}{2}-\alpha+\beta}}, \quad y \to 0, \Re(\beta)=0, \beta\ne 0
\end{align}
by \cite[9.220.2-4, 9.210.1]{GradRyzh} or \cite[13.14.18]{DLMF}.
Similarly, the functional equations \eqref{eq:Whitt1FEs}, which were computed by evaluating the Jacquet integral, can be verified from the functional equations \cite[Proposition 3.3]{HWI}.
For instance,
\begin{align*}
	-2 \Lambda_{(1,0,1)}(\mu) \bu^{1,-}_1 W^1(g, \mu, \psi_{1,1}) =& W^{1*}\paren{g, \mu^{w_4}} = \Lambda_{(0,1,1)}(\mu^{w_4}) \bu^{1,-}_0 W^1(g, \mu^{w_4}, \psi_{1,1}) \\
	=& \sqrt{2} \Lambda_{(0,1,1)}(\mu^{w_4}) \bu^{1,-}_0 T(w_5,\mu^{w_4}) W^1(g, \mu, \psi_{1,1}),
\end{align*}
so we conclude that $-2 \Lambda_{(1,0,1)}(\mu) \bu^{1,-}_1$ must be equal to $\sqrt{2} \Lambda_{(0,1,1)}(\mu^{w_4}) \bu^{1,-}_0 T^1(w_5,\mu^{w_4})$, which can be checked.

\subsection{The Kloosterman sums}
\label{sect:Kloos}
Define two exponential sums
\begin{align*}
	\tilde{S}(m_1,n_1,n_2;D_1,D_2) := \sum_{\substack{C_1 (\text{mod }D_1), C_2 (\text{mod }D_2)\\(C_1,D_1)=(C_2,D_2/D_1)=1}} e\left(n_1\frac{\bar{C_1}C_2}{D_1}+n_2\frac{\bar{C_2}}{D_2/D_1}+m_1\frac{C_1}{D_1}\right),
\end{align*}
for $D_1 | D_2$, and
\begin{align*}
	&S(m_1, m_2, n_1, n_2; D_1, D_2) \\
	& = \sum_{\substack{B_1, C_1 \summod{D_1}\\B_2, C_2 \summod{D_2}\\ }} \e{\frac{m_1B_1 + n_1(Y_1 D_2 - Z_1 B_2)}{D_1} + \frac{m_2B_2 + n_2(Y_2 D_1 - Z_2B_1)}{D_2}},
\end{align*}
where the sum is restricted to
\[ D_1C_2 + B_1B_2 + D_2C_1 \equiv 0 \summod{D_1D_2}, \qquad (B_1, C_1, D_1) = (B_2, C_2, D_2) = 1, \]
and the $Y_i$ and $Z_i$ are defined by
\[ Y_1B_1 + Z_1C_1 \equiv 1 \pmod{D_1}, \qquad Y_2B_2 + Z_2C_2 \equiv 1 \pmod{D_2}. \]
Then the $SL(3,\Z)$ Kloosterman sums are
\begin{align*}
	S_{w_4}(\psi_m,\psi_n;c) =& \delta_{\substack{n_2 c_1=m_1 c_2^2 \\ c_2|c_1}} \tilde{S}(-n_2,m_2,m_1;c_2,c_1), \\
	S_{w_5}(\psi_m,\psi_n;c) =& \delta_{\substack{n_1 c_2=m_2 c_1^2 \\ c_1|c_2}} \tilde{S}(n_1,m_1,m_2;c_1,c_2) \\
	S_{w_l}(\psi_m,\psi_n;c) =& S(-n_2,-n_1,m_1,m_2;c_1, c_2).
\end{align*}
Here $\delta_P$ is one if $P$ is true and zero otherwise.
Beware the previous note about the signs of the indices on $S_{w_l}$ (see the introduction to section \ref{sect:Background}).

We use the following bounds \cite{Stevens, BFG} (see \cite{Me01})
\begin{align}
\label{eq:Sw4Bound}
	\abs{S_{w_4}(\psi_m,\psi_n,c)} \le& d(c_1) (\abs{m_2}, \abs{n_2}, c_2) c_1, \\
	\abs{S_{w_5}(\psi_m,\psi_n,c)} \le& d(c_2) (\abs{m_2}, \abs{n_1}, c_1) c_2, \\
\label{eq:SwlBound}
	\abs{S_{w_l}(\psi_m,\psi_n,c)}^2 \le& d(c_1)^2 d(c_2)^2 \paren{\abs{m_1 n_2}, D} \paren{\abs{m_2 n_1}, D} (c_1, c_2) c_1 c_2,
\end{align}
where $D = \frac{c_1 c_2}{(c_1, c_2)}$.

\subsection{The spectral expansion}
For a Schwartz-class function $f:\Gamma\backslash G \to \C^3$ satisfying $f(gk)=f(g)\WigDMat{1}(k)$ which is orthogonal to the residual spectrum, the spectral expansion of \cite{HWI, HWII} takes the form
\begin{align}
\label{eq:InitSpectralExpand}
	f(g) =& \sum_{\varphi \in \mathcal{S}_3^1} \varphi(g) \, \int_{\Gamma\backslash G} f(g') \wbar{\trans{\varphi(g')}} dg' \\
		& \qquad +\frac{1}{2\pi i} \sum_{\phi \in \mathcal{S}_2^1} \sum_{\abs{m'}\le 1} \int_{\Re(\mu_1)=0} E^1_{m'}(g, \phi, \mu_1) \, \int_{\Gamma\backslash G} f(g') \wbar{\trans{E^1_{m'}(g', \phi, \mu_1)}} dg' d\mu_1 \nonumber \\
		& \qquad +\frac{1}{24 (2\pi i)^2} \sum_{\abs{m'}\le 1} \int_{\Re(\mu)=0} E^1_{m'}(g, \mu) \, \int_{\Gamma\backslash G} f(g') \wbar{\trans{E^1_{m'}(g', \mu)}} dg' d\mu, \nonumber
\end{align}
where $\mathcal{S}_3^1$ is a basis of vector-valued cusp forms for $\Gamma\backslash G$ and $\mathcal{S}_2^1$ is a basis of spherical Maass cusp forms for $SL(2,\Z)\backslash PSL(2,\R)$ (see the display before Theorem 1.1 in \cite{HWI}).
The Eisenstein series used here are defined in section 5 of \cite{HWI}.
Note that there is no need to consider the form $\wtilde{\varphi}$ described in \cite{HWII} since we are only dealing with the minimal $K$-type forms and the Whittaker function we employ is completed with respect to the $w_2$ functional equation.
Thus, in place of \cite[(166)]{HWII}, we are using the usual normalization of the cusp forms $\left<\varphi,\varphi\right>=1$.

The minimal parabolic Eisenstein series $E^1(g,\mu)$ (and the residual spectrum) and the maximal parabolic Eisenstein series $E^1(g,\phi,\mu_1)$ attached to an even $SL(2,\Z)$ Maass cusp form are both zero since their Fourier-Whittaker expansions involve the matrix $\Sigma^1_{++}=0$.
For the maximal parabolic Eisenstein series attached to an odd Maass cusp form, only the central row contributes, and we note that our choice of Whittaker functions is consistent, comparing \cite[(5.18),(5.19)]{HWI} to \eqref{eq:Whitt1FEs} since $\paren{\Sigma^1_{+-}}_0 = \bu^{1,-}_0$.

Therefore,
\begin{align}
\label{eq:SpectralExpand}
	f(g) =& \sum_{\varphi \in \mathcal{S}_3^1} \varphi(g) \int_{\Gamma\backslash G} f(g') \wbar{\trans{\varphi(g')}} dg' \\
		& \qquad +\frac{1}{2\pi i} \sum_{\substack{\phi \in \mathcal{S}_2^1\\ \phi\text{ odd}}} \int_{\Re(\mu_1)=0} E^1_0(g, \phi, \mu_1) \, \int_{\Gamma\backslash G} f(g') \wbar{\trans{E^1_0(g', \phi, \mu_1)}} dg' d\mu_1. \nonumber
\end{align}

\subsection{Integrals of gamma functions}
We will make use of Barnes' first and second lemmas:
\begin{thm}[Barnes' first lemma, {\cite[sect. 1.7]{Bailey}}]
\label{thm:BarnesFirst}
For $a,b,c,d\in\C$,
\begin{align*}
	\int_{-i\infty}^{+i\infty} \Gamma(a+s) \Gamma(b+s) \Gamma(c-s) \Gamma(d-s) \frac{ds}{2\pi i} &= \frac{\Gamma(a+c) \Gamma(b+c) \Gamma(a+d) \Gamma(b+d)}{\Gamma(a+b+c+d)}.
\end{align*}
\end{thm}
\begin{thm}[Barnes' second lemma, {\cite[sect. 6.2]{Bailey}}]
\label{thm:BarnesSecond}
For $a,b,c,d,e,f\in\C$ with $a+b+c+d+e-f=0$,
\begin{align*}
	& \int_{-i\infty}^{+i\infty} \frac{\Gamma(a+s) \Gamma(b+s) \Gamma(c+s) \Gamma(d-s) \Gamma(e-s)}{\Gamma(f+s)} \frac{ds}{2\pi i} \\
	&= \frac{\Gamma(a+d) \Gamma(b+d) \Gamma(c+d) \Gamma(a+e) \Gamma(b+e) \Gamma(c+e)}{\Gamma(f-a) \Gamma(f-b) \Gamma(f-c)}.
\end{align*}
\end{thm}

\subsection{The differential equations satisfied by the Kuznetsov kernel functions}
\label{sect:KernelDiffEqs}
We need to adjust the results of the paper \cite{SpectralKuz} to reflect the current notation.
That paper determines explicitly the power series solutions to the differential equations
\[ \Delta_i K_w(g,\mu) = \lambda_i(\mu) K_w(g,\mu), \qquad K_w(u g (wu'w^{-1}),\mu) = \psi_{1,1}(uu') K_w(g,\mu), \]
where $g\in G$, $w\in W$, $u\in U(\R)$, $u'\in U_w(\R)$ and
\[ \lambda_1(\mu) = 1-\tfrac{\mu_1^2+\mu_2^2+\mu_3^2}{2}, \qquad \lambda_2(\mu) = \mu_1 \mu_2 \mu_3. \]
We will only need the solutions in the specific cases $w=w_4,w_l$.

When $w=w_l$, the power-series solutions are
\[ J_{w_l}(y,\mu) = \abs{4\pi^2 y_1}^{1-\mu_3} \abs{4\pi^2 y_2}^{1+\mu_1} \sum_{n_1,n_2\ge 0} \frac{\Gamma\paren{n_1+n_2+\mu_1-\mu_3+1} \, (4\pi^2 y_1)^{n_1} (4\pi^2 y_2)^{n_2}}{\prod_{i=1}^3 \Gamma\paren{n_1+\mu_i-\mu_3+1}\Gamma\paren{n_2+\mu_1-\mu_i+1}}, \]
and these satisfy
\begin{align}
\label{eq:JwlAsymp}
	J_{w_l}(y,\mu) \sim p_{\rho+\mu}(y) \frac{(4\pi^2)^{2+\mu_1-\mu_3}}{\Gamma\paren{1+\mu_1-\mu_3} \Gamma\paren{1+\mu_1-\mu_2} \Gamma\paren{1+\mu_2-\mu_3}},
\end{align}
as $y\to 0$.

For $w=w_4$, we note that $\wtilde{K_w}(g,\mu):=K_w(g\vpmpm{+-},\mu)$ satisfies
\[ \wtilde{K_w}(g (wu'w^{-1}),\mu) = \psi_{-1,-1}(u') \wtilde{K_w}(g,\mu), \]
which matches the definition in \cite{SpectralKuz}, so we may use the solutions from that paper by replacing $y \mapsto y\vpmpm{+-}$:
\[ J_{w_4}(y,\mu) = \abs{8\pi^3 y_1}^{1-\mu_3} \sum_{n=0}^\infty \frac{(-8\pi^3 i y_1)^n}{n! \, \Gamma\paren{n+1+\mu_1-\mu_3} \, \Gamma\paren{n+1+\mu_2-\mu_3}}, \]
on the subspace of $U(\R) Y U_{w_4}(\R)$ (which is already smaller than $G$) defined by $y_2 = 1$.
These satisfy
\begin{align}
\label{eq:Jw4Asymp}
	J_{w_4}(y,\mu) \sim \abs{y_1}^{1-\mu_3} \frac{\abs{8\pi^3}^{1-\mu_3}}{\Gamma\paren{1+\mu_1-\mu_3} \Gamma\paren{1+\mu_2-\mu_3}},
\end{align}
as $y_1 \to 0$ with $y_2=1$.

\subsection{The Mellin-Barnes integrals for the Kuznetsov kernel functions}
\label{sect:OldMBIntegrals}

We will require the Mellin-Barnes integrals of these functions for the proof of Theorem \ref{thm:TechnicalWeyl}.

Define $W_3=\set{I,w_4,w_5}$, and
\begin{align}
\label{eq:d0SpectralMeasures}
	\cosmu(\mu) = \prod_{j<k} \cos \frac{\pi}{2}(\mu_j-\mu_k), \qquad \sinmu(\mu) = \prod_{j<k} \sin \frac{\pi}{2}(\mu_j-\mu_k).
\end{align}

If $y_1,y_2 > 0$,
\begin{align*}
	K_{w_l}(y,\mu) =& -\frac{\pi^3}{32} \sum_{w\in W} \frac{1}{\sinmu(\mu^w)} J_{w_l}(y,\mu^w) \\
	=& \frac{1}{4} \cosmu(\mu) \int_{\Re(s) = (2,2)} \wtilde{G}(0, 2s, 2\mu) (4\pi^2 y_1)^{1-s_1} (4\pi^2 y_2)^{1-s_2} \frac{ds}{(2\pi i)^2},
\end{align*}
where
\[ \wtilde{G}(0, s, \mu) := \frac{\Gamma\paren{\frac{s_1-\mu_1}{2}}\Gamma\paren{\frac{s_1-\mu_2}{2}}\Gamma\paren{\frac{s_1-\mu_3}{2}}\Gamma\paren{\frac{s_2+\mu_1}{2}}\Gamma\paren{\frac{s_2+\mu_2}{2}}\Gamma\paren{\frac{s_2+\mu_3}{2}}}{\Gamma\paren{\frac{s_1+s_2}{2}}}. \]

If $y_1,y_2 < 0$,
\begin{align*}
	K_{w_l}^{--}(y, \mu) :=& J_{w_l}(y,\mu) - J_{w_l}(y,\mu^{w_l}) \\
	=& -\frac{1}{\pi} \sin\pi(\mu_1-\mu_3) \int_{-i\infty}^{+i\infty} \int_{-i\infty}^{+i\infty} \abs{4\pi^2 y_1}^{1-s_1} \abs{4\pi^2 y_2}^{1-s_2} \\
	& \qquad \times \frac{\Gamma\paren{s_1-\mu_3} \Gamma\paren{s_1-\mu_1} \Gamma\paren{s_2+\mu_1} \Gamma\paren{s_2+\mu_3}}{\Gamma\paren{1-s_1+\mu_2} \Gamma\paren{s_1+s_2} \Gamma\paren{1-s_2-\mu_2}} \frac{ds_1}{2\pi i} \frac{ds_2}{2\pi i}.
\end{align*}

If $y_1 < 0 < y_2$,
\begin{align*}
	K_{w_l}^{-+}(y, \mu) :=& J_{w_l}(y,\mu)-J_{w_l}(y,\mu^{w_2}) \\
	=& -\frac{1}{\pi} \sin\pi(\mu_1-\mu_2) \int_{-i\infty}^{+i\infty} \int_{-i\infty}^{+i\infty} \abs{4\pi^2 y_1}^{1-s_1} \abs{4\pi^2 y_2}^{1-s_2} \\
	& \qquad \times \frac{\Gamma\paren{s_1-\mu_3} \Gamma\paren{s_2+\mu_1} \Gamma\paren{s_2+\mu_2} \Gamma\paren{1-s_2-s_1}}{\Gamma\paren{1+\mu_1-s_1} \Gamma\paren{1+\mu_2-s_1} \Gamma\paren{1-s_2-\mu_3}} \frac{ds_1}{2\pi i} \frac{ds_2}{2\pi i}.
\end{align*}

If $y_1 > 0 > y_2$,
\begin{align*}
	K_{w_l}^{+-}(y, \mu) :=& J_{w_l}(y,\mu) - J_{w_l}(y,\mu^{w_3}) \\
	=& -\frac{1}{\pi} \sin\pi(\mu_2-\mu_3) \int_{-i\infty}^{+i\infty} \int_{-i\infty}^{+i\infty} \abs{4\pi^2 y_1}^{1-s_1} \abs{4\pi^2 y_2}^{1-s_2} \\
	& \qquad \times \frac{\Gamma\paren{s_1-\mu_2} \Gamma\paren{s_1-\mu_3}\Gamma\paren{s_2+\mu_1} \Gamma\paren{1-s_1-s_2}}{\Gamma\paren{1-s_1+\mu_1}\Gamma\paren{1-\mu_2-s_2} \Gamma\paren{1-\mu_3-s_2}} \frac{ds_2}{2\pi i} \frac{ds_1}{2\pi i}.
\end{align*}

\section{Stade's Formula}
We prove Theorem \ref{thm:StadesFormula}.
By parity considerations, we need only consider the terms
\begin{align*}
	16\pi^{3t} \Psi^1 =& 2\Psi^1_{00}+2\Psi^1_{10}+2\Psi^1_{11}, \\
	\Psi^1_{\ell_1,\ell_2} :=& \int_{\Re(s)=\mathfrak{s}} \wtilde{G}(1,(1-\ell_1,1-\ell_1,\ell_1),(\ell_2,\ell_2,1-\ell_2),s,\mu) \\
	& \qquad \times \wtilde{G}(1,(1-\ell_1,1-\ell_1,\ell_1),(\ell_2,\ell_2,1-\ell_2), (2t-s_1,t-s_2),\mu') \frac{ds}{(2\pi i)^2}.
\end{align*}

We insert the definition \eqref{eq:GDef}, and collect terms as follows:
\begin{align*}
	\Psi^1_{\ell_1,\ell_2} =& \int_{\Re(s)=\mathfrak{s}} \frac{\Gamma\paren{\tfrac{1-\ell_1+s_1-\mu_1}{2}} \Gamma\paren{\tfrac{1-\ell_1+s_1-\mu_2}{2}} \Gamma\paren{\tfrac{\ell_2+s_2+\mu_1}{2}} \Gamma\paren{\tfrac{\ell_2+s_2+\mu_2}{2}}}{\Gamma\paren{\tfrac{1-\ell_1+\ell_2+s_1+s_2}{2}}} \\
	& \qquad \times \frac{\Gamma\paren{\tfrac{1-\ell_1+2t-s_1-\mu_1'}{2}} \Gamma\paren{\tfrac{1-\ell_1+2t-s_1-\mu_2'}{2}} \Gamma\paren{\tfrac{\ell_2+t-s_2+\mu_1'}{2}} \Gamma\paren{\tfrac{\ell_2+t-s_2+\mu_2'}{2}}}{\Gamma\paren{\tfrac{1-\ell_1+\ell_2+3t-s_1-s_2}{2}}} \\
	& \qquad \times \Gamma\paren{\tfrac{\ell_1+s_1-\mu_3}{2}} \Gamma\paren{\tfrac{1-\ell_2+s_2+\mu_3}{2}} \Gamma\paren{\tfrac{\ell_1+2t-s_1-\mu_3'}{2}} \Gamma\paren{\tfrac{1-\ell_2+t-s_2+\mu_3'}{2}} \frac{ds}{(2\pi i)^2}.
\end{align*}

To both of the quotients in the integrand, we apply Barnes' first lemma in reverse, giving
\begin{align*}
	\Psi^1_{\ell_1,\ell_2} =& \frac{1}{4} \int_{\Re(u)=\mathfrak{u}} \int_{\Re(s)=\mathfrak{s}} \Gamma\paren{\tfrac{1-\ell_1+s_1+\mu_3+u_1}{2}} \Gamma\paren{\tfrac{\ell_2+s_2+u_1}{2}} \Gamma\paren{\tfrac{\mu_1-u_1}{2}} \Gamma\paren{\tfrac{\mu_2-u_1}{2}} \\
	& \qquad \times \Gamma\paren{\tfrac{1-\ell_1+t-s_1+\mu_3'+u_2}{2}} \Gamma\paren{\tfrac{\ell_2-s_2+u_2}{2}} \Gamma\paren{\tfrac{t+\mu_1'-u_2}{2}} \Gamma\paren{\tfrac{t+\mu_2'-u_2}{2}} \\
	& \qquad \times \Gamma\paren{\tfrac{\ell_1+s_1-\mu_3}{2}} \Gamma\paren{\tfrac{1-\ell_2+s_2+\mu_3}{2}} \Gamma\paren{\tfrac{\ell_1+2t-s_1-\mu_3'}{2}} \Gamma\paren{\tfrac{1-\ell_2+t-s_2+\mu_3'}{2}} \frac{ds}{(2\pi i)^2} \frac{du}{(2\pi i)^2}.
\end{align*}

Now apply Barnes' first lemma to the $s$-integrals,
\begin{align*}
	\Psi^1_{\ell_1,\ell_2} =& \int_{\Re(u)=\mathfrak{u}} \frac{\Gamma\paren{\tfrac{2-2\ell_1+t+\mu_3+\mu_3'+u_1+u_2}{2}} \Gamma\paren{\tfrac{1+2t+\mu_3-\mu_3'+u_1}{2}} \Gamma\paren{\tfrac{1+t-\mu_3+\mu_3'+u_2}{2}} \Gamma\paren{\tfrac{2\ell_1+2t-\mu_3-\mu_3'}{2}}}{\Gamma\paren{\tfrac{2+3t+u_1+u_2}{2}}} \\
	& \qquad \times \frac{\Gamma\paren{\tfrac{2\ell_2+u_1+u_2}{2}} \Gamma\paren{\tfrac{1+t+\mu_3'+u_1}{2}} \Gamma\paren{\tfrac{1+\mu_3+u_2}{2}} \Gamma\paren{\tfrac{2-2\ell_2+t+\mu_3+\mu_3'}{2}}}{\Gamma\paren{\tfrac{2+t+\mu_3+\mu_3'+u_1+u_2}{2}}} \\
	& \qquad \times \Gamma\paren{\tfrac{\mu_1-u_1}{2}} \Gamma\paren{\tfrac{\mu_2-u_1}{2}} \Gamma\paren{\tfrac{t+\mu_1'-u_2}{2}} \Gamma\paren{\tfrac{t+\mu_2'-u_2}{2}} \frac{ds}{(2\pi i)^2} \frac{du}{(2\pi i)^2}.
\end{align*}

When $\ell_1=\ell_2=0$, we may cancel the factors $\Gamma\paren{\frac{2+t+\mu_3+\mu_3'+u_1+u_2}{2}}$, but for the remaining two terms we use
\[ \Gamma\paren{\tfrac{u_1+u_2}{2}} \Gamma\paren{\tfrac{2+t+\mu_3+\mu_3'}{2}}+ \Gamma\paren{\tfrac{2+u_1+u_2}{2}} \Gamma\paren{\tfrac{t+\mu_3+\mu_3'}{2}} = \tfrac{t+\mu_3+\mu_3'+u_1+u_2}{2}\Gamma\paren{\tfrac{u_1+u_2}{2}} \Gamma\paren{\tfrac{t+\mu_3+\mu_3'}{2}}, \]
and so the same factor cancels in the sum.
Similarly, we have
\[ \Gamma\paren{\tfrac{2t-\mu_3-\mu_3'}{2}} \Gamma\paren{\tfrac{2+t+\mu_3+\mu_3'}{2}} + \Gamma\paren{\tfrac{2+2t-\mu_3-\mu_3'}{2}} \Gamma\paren{\tfrac{t+\mu_3+\mu_3'}{2}} = \tfrac{3t}{2} \Gamma\paren{\tfrac{2t-\mu_3-\mu_3'}{2}} \Gamma\paren{\tfrac{t+\mu_3+\mu_3'}{2}}, \]
and we are left evaluate
\begin{align*}
	16\pi^{3t} \Psi^1 =& 3t \Gamma\paren{\tfrac{2t-\mu_3-\mu_3'}{2}} \Gamma\paren{\tfrac{t+\mu_3+\mu_3'}{2}} \\
	& \qquad \times \int_{\Re(u)=\mathfrak{u}} \frac{\Gamma\paren{\tfrac{1+2t+\mu_3-\mu_3'+u_1}{2}} \Gamma\paren{\tfrac{u_1+u_2}{2}} \Gamma\paren{\tfrac{1+t+\mu_3'+u_1}{2}} \Gamma\paren{\tfrac{\mu_1-u_1}{2}} \Gamma\paren{\tfrac{\mu_2-u_1}{2}}}{\Gamma\paren{\tfrac{2+3t+u_1+u_2}{2}}} \\
	& \qquad \times \Gamma\paren{\tfrac{1+\mu_3+u_2}{2}} \Gamma\paren{\tfrac{1+t-\mu_3+\mu_3'+u_2}{2}} \Gamma\paren{\tfrac{t+\mu_1'-u_2}{2}} \Gamma\paren{\tfrac{t+\mu_2'-u_2}{2}} \frac{ds}{(2\pi i)^2} \frac{du}{(2\pi i)^2}.
\end{align*}
The result follows from two applications of Barnes' second lemma.

\section{Kontorovich-Lebedev Inversion}
If $F(\mu)$ is holomorphic in a neighborhood $\abs{\Re(\mu)} < \delta < \frac{1}{10}$, then we can argue as in \cite{GoldKont} that the $Y^+$ integral of $(F^\flat)^\sharp$ converges absolutely and define
\[ F(\mu,\epsilon) := \int_{Y^+} F^\flat(y) \wbar{\trans{W^{1*}(y,\mu)}} (y_1^2 y_2)^\epsilon dy = \int_{\Re(\mu')=0} F(\mu') \Psi^1(\mu',-\mu,\epsilon) \sinmu^1(\mu')\,d\mu', \]
where we assume $0<\epsilon<\frac{\delta}{2} =: \eta$, and $\Re(\mu)=0$ with $\mu_1\ne\mu_2$.
The points with $\mu_1=\mu_2$ and $\Re(\mu)\ne0$ may be handled by analytic continuation.
There is a slight subtlety that $\wbar{\mu}$ is either $-\mu$ or $-\mu^{w_2}$, but $\Psi^1$ is invariant under $\mu' \mapsto (\mu')^{w_2}$.

Shift the $\mu'$ integral to $\Re(\mu')=(-\eta,-\eta,2\eta)$, giving
\begin{align*}
	F(\mu,\epsilon) =& \int_{\Re(\mu')=(-\eta,-\eta,2\eta)} F(\mu') \Psi^1(\mu',-\mu,\epsilon) \sinmu^1(\mu')\,d\mu' \\
	& + (2\pi i)^2 \res_{\mu_1'=\mu_1-\epsilon} \res_{\mu_2'=\mu_2-\epsilon} F(\mu') \sinmu^1(\mu') \Psi^1(\mu',-\mu,\epsilon) \\
	& + (2\pi i)^2 \res_{\mu_1'=\mu_2-\epsilon} \res_{\mu_2'=\mu_1-\epsilon} F(\mu') \sinmu^1(\mu') \Psi^1(\mu',-\mu,\epsilon) \\
	& + \text{ mixed terms}.
\end{align*}
Note that $\sinmu^1(\mu')=0$ when $\mu_1'=\mu_2'$.
Evaluating the residues and taking the limit $\epsilon \to 0$ gives the result.

\section{Kuznetsov's Formula}
For convenience, we abbreviate the spectral expansion \eqref{eq:SpectralExpand} in the form
\[ f(g) = \int_\mathcal{B} \xi(g) \, \int_{\Gamma\backslash G} f(g') \wbar{\trans{\xi(g')}} dg' \, d\xi, \]
where $\mathcal{B}$ stands for the combined cuspidal and continuous spectral bases.
Applying the expansion to \eqref{eq:MainPoincare} with Theorem \ref{thm:KontLebedev}, we have
\begin{align}
\label{eq:PreKuzSpectral}
	\mathcal{P}(yk) :=& \int_{U(\Z)\backslash U(\R)} P_m(xyk,F) \wbar{\psi_n(x)} dx \\
	=& \abs{\frac{m_1 m_2}{n_1 n_2}} \int_\mathcal{B} F(\mu_\xi) \wbar{\rho_\xi^*(m)} \rho_\xi^*(n) W^{1*}(\wtilde{n} yk, \mu_\xi) \, d\xi, \nonumber
\end{align}
where we define the Fourier-Whittaker coefficients of a Maass form $\xi$ with Langlands parameters $\mu_\xi$ by
\begin{align}
\label{eq:FWcoefDef}
	\int_{U(\Z)\backslash U(\R)} \xi(xyk) \wbar{\psi_m(x)} dx = \frac{\rho_\xi^*(n)}{\abs{m_1 m_2}} W^{1*}(\wtilde{m} yk, \mu_\xi).
\end{align}
Our use of the completed Whittaker function here implies the Fourier-Whittaker coefficients $\rho_\xi^*(m)$ contain an exponential growth factor, which will be offset later using Stade's formula.
We are assuming $m_1 m_2 n_1 n_2 \ne 0$ for simplicity.

Define the integral transform
\begin{align}
\label{eq:HwTildeDef}
	\wtilde{H}_w(F; y, g) =& \frac{1}{\abs{y_1 y_2}} \int_{\wbar{U}_w(\R)} \int_{\Re(\mu)=0} F(\mu) W^{1*}(ywxg,\mu) \sinmu^1(\mu)\,d\mu \,\wbar{\psi_{1,1}(x)} dx,
\end{align}
for $w\in W$, $y\in Y:=VY^+ \cong (\R^\times)^2$, $g\in G$.
Then the Bruhat decomposition (in the form $U(\Q)CWU(\Q)V$) applied to the sum over $U(\Z)\backslash\Gamma$ in \eqref{eq:MainPoincare} implies (see \cite[section 10.6]{Gold01} or \cite{BFG} and the substitutions in \cite[section 4.2]{Me01})
\begin{align}
\label{eq:PreKuzBruhat}
	\mathcal{P}(yk) =& \sum_{w\in W} \sum_{v\in V} \sum_{c_1,c_2\ge1} \frac{S_w(\psi_m,\psi_n^v,c)}{c_1 c_2} \abs{\frac{m_1 m_2}{n_1 n_2}} \wtilde{H}_w\paren{F; \wtilde{m} c w v \wtilde{n}^{-1} w^{-1},\wtilde{n} yk}.
\end{align}
The equality of \eqref{eq:PreKuzSpectral} and \eqref{eq:PreKuzBruhat} is called the pre-Kuznetsov formula.

Note: To be precise, in the development of the Kuznetsov formula, we must initially require $F$ to be holomorphic on $\abs{\Re(\mu_i)} < 1+\delta$ with zeros also at $\mu_1-\mu_3=0,\pm2$ and $\mu_2-\mu_3=0,\pm2$ and exponential decay to overcome the growth of the Fourier-Whittaker coefficients for absolute convergence of the Poincar\'e series (see the proof of lemma \ref{lem:KuzKernels} below).
We may relax to $\abs{\Re(\mu_i)} < \frac{1}{2}+\delta$ once we reach the pre-Kuznetsov formula, using the bound \eqref{eq:SwlBound}, and the extra zeros and exponential decay will come from the final application of Stade's formula.

In section \ref{sect:KuzKernelsProof} below, we show
\begin{lem}
\label{lem:KuzKernels}
	Let $F$ be holomorphic and Schwartz-class on a neigborhood of $\Re(\mu)=0$, then for $w=I,w_4,w_5,w_l$, we have
	\begin{align*}
		\wtilde{H}_w(F; y, g) =& \frac{1}{\abs{y_1 y_2}} \int_{\Re(\mu)=0} F(\mu) K_w^1(y,\mu) W^{1*}(g,\mu) \sinmu^1(\mu)\,d\mu,
	\end{align*}
	with $K_w^1(y,\mu)$ as in \eqref{eq:K1I}-\eqref{eq:K1wl}.
\end{lem}
The terms corresponding to the Weyl elements $w_2$ and $w_3$ disappear by the compatibility condition and the assumption of non-degenerate characters (see the discussion in \cite[section 1.7]{SpectralKuz}), but a similar statement holds in those cases, as well.

Then we integrate away the remaining Whittaker function by inserting a factor
\[ W^{1*}_0(y'k',-\mu) (y_1')^{2t} (y_2')^t \]
into \eqref{eq:MainPoincare} (that is, substituting on $F(\mu)$), evaluating the pre-Kuznetsov formula at $yk=\wtilde{n}^{-1}y'k'$, taking the central entry, integrating over $y'k'$ using Stade's formula at $t>1$ and applying dominated convergence to get to $t=1$.
The appearance of $\cosmu^1(\mu)$ in this particular method of deriving the Kuznetsov formula may seem somewhat artificial, but it will be clear this is the correct weighting when we move to bases of Hecke eigenfunctions in Theorem \ref{thm:HeckeKuznetsov}.

For the Eisenstein series terms, recall \cite[(5.6)]{HWI}; even though $\mu_\phi$ (for $\phi \in \mathcal{S}^1_2$) is only defined up to sign, both $F(\mu)$ and $\cosmu^1(\mu)$ are invariant under $\mu \mapsto \mu^{w_2}$.

\subsection{The proof of Lemma \ref{lem:KuzKernels}}
\label{sect:KuzKernelsProof}
Note that no work need be done for the trivial Weyl element, so we need only determine the kernel function for $w=w_4,w_5,w_l$.
The $y$ argument of $\wtilde{H}_w$ is necessarily signed, so we replace it with $y\vpmpm{\varepsilon_1,\varepsilon_2}$ where now $y\in Y^+$ and $\vpmpm{\varepsilon_1,\varepsilon_2}\in V$.

Write $W_3=\set{I,w_4,w_5}$, as before.
With $c_W(w)$, $\alpha(w)$ and $\bu(w)$ as in section \ref{sect:WhittFuncs}, let $\wtilde{\alpha}(w)=\alpha(w)^{w_l}-(1,1,1)$.
Then assuming the entries of $\mu$ are distinct, we may write
\begin{align}
\label{eq:WhittResExpand}
	W^{1*}(y, \mu) =& \pi^{\frac{3}{2}} \sum_{w\in W} \abs{c_W(w)} p_{\rho+\mu^w}(y) \Lambda_{\wtilde{\alpha}(w)}(\mu^{w w_l}) \bu(w) \\
	&+\frac{1}{4\pi^2} \sum_{w\in W_3} (\pi y_1)^{1-\mu_3^w} \int_{\Re(s_2)=-\eta}  (\pi y_2)^{1-s_2} \res_{s_1=\mu_3^w} G^1\paren{s,\mu} \frac{ds_2}{2\pi i} \nonumber\\
	&+\frac{1}{4\pi^2} \sum_{w\in W_3} (\pi y_2)^{1+\mu_1^w} \int_{\Re(s_1)=-\eta} (\pi y_1)^{1-s_1} \res_{s_2=-\mu_1^w} G^1\paren{s,\mu} \frac{ds_1}{2\pi i} \nonumber\\
	&+\frac{1}{4\pi^2} \int_{\Re(s)=(-\eta,-\eta)} (\pi y_1)^{1-s_1} (\pi y_2)^{1-s_2} G^1\paren{s,\mu} \frac{ds}{(2\pi i)^2}, \nonumber
\end{align}
for some small $\eta > 0$.
The particular form of the secondary terms is not essential here.

In \eqref{eq:MainPoincare}, we may apply the above expansion of the Whittaker function, and then shift the contour in the $\mu$ integral term-by-term until the real part of the powers of $y_1$ and $y_2$ are nearly equal and greater than 1; this guarantees the absolute convergence of the $x$ integral in \eqref{eq:HwTildeDef}.
The larger region of holomorphy and zeros for $F(\mu)$ in the Kuznetsov formula similarly guarantee the absolute convergence of the sum of Kloosterman sums (see the note after \eqref{eq:PreKuzBruhat}).

By way of refining the process, we drop the $K^\text{as}_w$ construction of \cite{SpectralKuz}.
For each term of \eqref{eq:WhittResExpand}, we may pull the $x$ integral of \eqref{eq:HwTildeDef} inside the $\mu$ and any possible $s$ integrals.
Up to permutations of $\mu$, and substituting coordinates of $s$ for coordinates of $\mu$ as necessary in the secondary terms, the isolated $x$ integral takes the form
\begin{align}
\label{eq:XDef}
	X:=& p_{\rho+\mu}(y) p_{\rho+\mu^w}(t) \WigDMat{1}(\vpmpm{\varepsilon_1,\varepsilon_2}) \int_{\wbar{U}_w(\R)} p_{\rho+\mu}(y^*) \WigDMat{1}(k^*) \psi_{yt^w}(x^*) \wbar{\psi_t(x)} dx,
\end{align}
for $g=t\in Y^+$, where $x^* y^* k^*=wx$.
The key technical difficulty is the analytic continuation of the integral $X$ to a neighborhood of $\Re(\mu)=0$.
We summarize the modifications to the argument of \cite[section 2.6.2]{SpectralKuz} in the case $w=w_l$ in section \ref{sect:LongEleXCont} below.

Temporarily replacing $F(\mu)$ with an approximation to the identity, say
\[ f(\mu)=f(\mu,\mu',\delta)=\delta^{-2} \exp\paren{\delta^{-2}\sum_{i=1}^3(\mu_i^w-\mu_i')^2}, \]
as in \cite[section 2.6.1]{SpectralKuz}, we see that
\begin{align}
\label{eq:tildeKwDef}
	\wtilde{H}_w(F;g,g') = \frac{1}{\abs{y_1 y_2}} \int_{\Re(\mu)=0} F(\mu) \wtilde{K}_w(g,g';\mu) \sinmu^1(\mu)\,d\mu,
\end{align}
for some vector-valued function $\wtilde{K}_w$ which satisfies
\begin{align}
\label{eq:tildeKwChary}
	\wtilde{K}_w(xg(w x' w),\cdot;\mu) = \psi_{1,1}(xx') \wtilde{K}_w(g,\cdot;\mu), \\
\label{eq:tildeKwChart}
	\wtilde{K}_w(\cdot,xg'k;\mu) = \psi_{1,1}(x) \wtilde{K}_w(\cdot, g';\mu)\WigDMat{1}(k),
\end{align}
and is an eigenfunction of both Casimir operators with eigenvalues matching $p_{\rho+\mu}$ in both $g$ and $g'$, for $g$ in the appropriate subspace of $G$ (see section \ref{sect:KernelDiffEqs}).
This subspace is determined by the requirement that \eqref{eq:HwTildeDef} satisfies the $x$-invariance condition of \eqref{eq:tildeKwChart}; e.g. $y_1=1$ for $w=w_4$.

In sections \ref{sect:LongEleAsymp}-\ref{sect:w5Asymp}, we compute the asymptotics of $\wtilde{K}_w(y,t;\mu)$ as one or both of the $y_i \to 0$.
In each case, the conclusion is
\[ \wtilde{K}_w(y\vpmpm{\varepsilon_1,\varepsilon_2},t;\mu) \sim K_w^1(y\vpmpm{\varepsilon_1,\varepsilon_2},\mu) W^{1*}(t,\mu) \]
as the relevant coordinates of $y$ tend to 0, where $K_w^1(y,\mu)$ solves the differential equations coming from the Casimir operators and \eqref{eq:tildeKwChary}.
Having already solved the requisite differential equations in section \ref{sect:KernelDiffEqs}, this is sufficient to identify the particular linear combination of power series solutions $J_w(y,\mu)$ by comparing the first-term asymptotics as $y\to0$, and we conclude that actually
\[ \wtilde{K}_w(y\vpmpm{\varepsilon_1,\varepsilon_2},t;\mu) = K_w^1(y\vpmpm{\varepsilon_1,\varepsilon_2},\mu) W^{1*}(t,\mu). \]

\subsubsection{Analytic continuation in the long element case}
\label{sect:LongEleXCont}
From the definition of the incomplete Whittaker function \cite[(3.3)]{HWI},
\begin{align}
\label{eq:XAsympExpand}
	X=& p_{\rho+\mu}(y) \WigDMat{1}(\vpmpm{\varepsilon_1,\varepsilon_2}) W^1(t,w,\mu,\psi_{1,1})+p_{\rho+\mu}(y) p_{\rho+\mu^w}(t) \WigDMat{1}(\vpmpm{\varepsilon_1,\varepsilon_2}) X_1, \\
	X_1 :=& \int_{\wbar{U}_w(\R)} p_{\rho+\mu}(y^*) \WigDMat{1}(k^*) (\psi_{yt^w}(x^*)-1) \wbar{\psi_t(x)} dx. \nonumber
\end{align}

We need to know that $X_1$ continues to a function which is sufficiently differentiable in $y$ and $t$ and asymptotically smaller than $1$ as $y \to 0$ (i.e. $O(y_1^\epsilon+y_2^\epsilon)$), but we will never directly use the integral representation by which we prove this.
If we can arrange that the $x$ integral converges very rapidly, the asymptotic condition on $y$ follows from the mean value theorem applied to $\psi_y(x^*)-1$.
Similarly, rapid convergence in $x$ means the differentiablity is satisfied as well.

As in \cite[section 2.2]{HWI}, let $\tildek{u_1,u_2,u_3}=k^*$.
Applying the trick \cite[(2.7)-(2.8)]{HWI}, we may write each component of $X_1$ as a finite linear combination of integrals of the form
\begin{align*}
	X_1' :=& \int_{\wbar{U}_w(\R)} p_{\rho+\mu}(y^*) u_1^{m_1} u_2^{m_2} u_3^{m_3} (\psi_{yt^w}(x^*)-1) \wbar{\psi_t(x)} dx,
\end{align*}
for some $m\in\Z^3$.
Out of a fervent desire to not repeat this part of the argument in the second paper, we point out the same is true if $\WigDMat{1}$ is replaced with $\WigDMat{d}$ in $X_1$.

The precise form of the $u_i$ may be found in \cite[(3.14)]{HWI}, but we only need note that each may be written as a finite sum of products of the terms
\[ x_1, \; x_2, \; x_3, \; \sqrt{1+x_2^2}, \; \sqrt{1+x_2^2+x_3^2}, \; \sqrt{1+x_1^2+(x_3 - x_1 x_2)^2}, \]
thus the argument of \cite[section 2.6.2]{SpectralKuz} goes through essentially unchanged, but with potentially many more iterations.
(A more careful analysis would reveal that no additional iterations of the integration-by-parts procedure are necessary as the derivatives of the $u_i$ generally have no singularities, even at infinity.)
The point is that the $X_3'$ integral considered there is sufficiently general to handle the current case, as well.

\subsubsection{Asymptotics in the long element case}
\label{sect:LongEleAsymp}
Set $W_2=\set{I,w_2}$.
Applying \eqref{eq:XAsympExpand} and \eqref{eq:WhittResExpand} to \eqref{eq:HwTildeDef} along with \eqref{eq:uSigns} and \eqref{eq:Whitt1FEs} implies
\begin{align*}
	\wtilde{K}_{w_l}(y,t;\mu) \sim \pi^{\frac{3}{2}} W^{1*}(t,\mu) \biggl(&\varepsilon_2 \sum_{w\in W_2} p_{\rho+\mu^w}(y) \frac{\Lambda_{(0,0,-1)}(\mu^{w w_l})}{\Lambda_{(0,1,1)}(\mu^w)} \\
	&+\varepsilon_1 \sum_{w\in W_2 w_4} p_{\rho+\mu^w}(y) \frac{\Lambda_{(-1,0,0)}(\mu^{w w_l})}{\Lambda_{(1,1,0)}(\mu^w)} \\
	&-\varepsilon_1 \varepsilon_2 \sum_{w\in W_2 w_5} p_{\rho+\mu^w}(y) \frac{\Lambda_{(0,-1,0)}(\mu^{w w_l})}{\Lambda_{(1,0,1)}(\mu^w)}\biggr),
\end{align*}
as $y \to 0$, and this matches $K_{w_l}^1(y\vpmpm{\varepsilon_1,\varepsilon_2},\mu) W^{1*}(t,\mu)$ as in \eqref{eq:K1wl} (compare \eqref{eq:JwlAsymp}).

\subsubsection{Asymptotics in the $w_4$ case}
\label{sect:w4Asymp}
By analytic continuation, we may replace \eqref{eq:WhittResExpand} with \eqref{eq:Whitty1Asymp} in the construction of $\wtilde{K}_{w_4}(y,t; \mu)$.
Inserting \eqref{eq:Whitty1Asymp} into \eqref{eq:HwTildeDef} and applying \cite[(3.11)]{HWI}, we see that
\begin{align*}
	& \wtilde{K}_{w_4}(y,t; \mu) \sim \\
	& \sum_{w\in W_3} c_W(w) \Lambda_{\alpha(w)}(\mu^w) y_1^{1-\mu_1^w} p_{\rho+\mu^{w w_3}}(t) \bu(w) \mathcal{W}^1(0,\mu_1^w-\mu_2^w) \WigDMat{1}(w_3) \mathcal{W}^1(0,\mu_1^w-\mu_3^w) \\
	& \qquad \times \WigDMat{1}(w_5 \vpmpm{\varepsilon_1,-1}) \int_{\wbar{U}_{w_4}(\R)} p_{\rho+\mu^{w w_l}}(y^*) \mathcal{W}^1(t_1 y_2^*,\mu_2^w-\mu_3^w) \WigDMat{1}(\vpmpm{+-} k^*) \psi_{0,t_1}(x^*) \wbar{\psi_t(x)} dx,
\end{align*}
as $y_1 \to 0$ on the subspace $y_2=\varepsilon_2=1$, where $w_4 x = x^* y^* k^*$.

We may compute (see \cite[section 2.4]{HWI})
\[ x_2^* = -\frac{x_2 x_3}{1+x_2^2}, \qquad y_1^*=\frac{\sqrt{1+x_2^2}}{1+x_2^2+x_3^2}, \qquad y_2^*=\frac{\sqrt{1+x_2^2+x_3^2}}{1+x_2^2}, \]
\[ \vpmpm{+-} k^*= \wtilde{k}\paren{-i,\tfrac{-x_3+i\sqrt{1+x_2^2}}{\sqrt{1+x_2^2+x_3^2}},\tfrac{1-ix_2}{\sqrt{1+x_2^2}}}. \]
Then, recognizing \cite[(3.22)]{HWI} and using \eqref{eq:uSigns} and \cite[(2.18)]{HWI}, we have
\begin{align*}
	\wtilde{K}_{w_4}(y,t; \mu) \sim&\pi W^{1*}(t,\mu) \biggl(i\varepsilon_1 y_1^{1-\mu_1} \frac{\Lambda_{(-1,0,1)}(\mu)}{\Lambda_{(1,0,1)}(\mu^{w_5})} -y_1^{1-\mu_3} \frac{\Lambda_{(0,0,0)}(\mu^{w_4})}{\Lambda_{(0,1,1)}(\mu)} +i\varepsilon_1 y_1^{1-\mu_2} \frac{\Lambda_{(0,-1,1)}(\mu^{w_5})}{\Lambda_{(1,1,0)}(\mu^{w_4})}\biggr),
\end{align*}
as $y_1 \to 0$ with $y_2=\varepsilon_2=1$, and this matches $K_{w_4}^1(y\vpmpm{\varepsilon_1,\varepsilon_2},\mu) W^{1*}(t,\mu)$ as in \eqref{eq:K1w4} (compare \eqref{eq:Jw4Asymp}).

\subsubsection{The $w_5$ case}
\label{sect:w5Asymp}
From \cite[section 6.3]{HWII} and the direct computation $\WigDRow{d}{0}(\vpmpm{--}w_l) = -\sqrt{2} \bu^{1,+}_1$, we have
\begin{align}
\label{eq:DualWhitt}
	W^{1*}(g, \mu) = -W^{1*}(\vpmpm{--}g^\iota w_l, -\mu),
\end{align}
where $g^\iota = w_l \trans{(g^{-1})} w_l$.
Then, as in \cite[section 2.3]{SpectralKuz}, \eqref{eq:HwTildeDef} becomes
\begin{align*}
	\wtilde{H}_{w_5}(F; y, g) =& -\frac{1}{\abs{y_1 y_2}} \int_{\wbar{U}_{w_4}(\R)} \int_{\Re(\mu)=0} F(\mu) \\
	& \times W^{1*}(\vpmpm{-+} y^\iota w_4x \vpmpm{--}g^\iota w_l,-\mu) \sinmu^1(\mu)\,d\mu \wbar{\psi_{1,1}(x)} dx,
\end{align*}
after sending $x^\iota \mapsto \vpmpm{--} x \vpmpm{--}$.
That is,
\begin{align*}
	\wtilde{H}_{w_5}(F; y, g) =& \wtilde{H}_{w_4}(\wtilde{F}; \vpmpm{-+} y^\iota, \vpmpm{--}g^\iota w_l), \\
	\wtilde{F}(\mu) =& -F(-\mu).
\end{align*}

From the conclusion of the lemma in the case $w=w_4$, we have
\begin{align*}
	\wtilde{H}_{w_5}(F; y, g) =& \frac{1}{\abs{y_1 y_2}} \int_{\Re(\mu)=0} \wtilde{F}(\mu) K^1_{w_4}(\vpmpm{-+} y^\iota;\mu) W^{1*}(\vpmpm{--}g^\iota w_l,\mu) \sinmu^1(\mu)\,d\mu \\
	=& \frac{1}{\abs{y_1 y_2}} \int_{\Re(\mu)=0} F(\mu) K^1_{w_4}(\vpmpm{-+} y^\iota,-\mu) W^{1*}(g,\mu) \sinmu^1(\mu)\,d\mu.
\end{align*}

\section{The Weyl Law}
\subsection{The proof of Theorem \ref{thm:TechnicalWeyl}}
We use the following lemma:
\begin{lem}
\label{lem:ElemBounds}
	For $v\in\R^3$ with $v_3>v_2>v_1$ and $v_3-v_2>v_2-v_1$,
	\[ \int_{\R^2} \abs{1+i(u_1-u_2)}^{\frac{3}{2}+\epsilon} \prod_{i=1}^3 \prod_{j=1}^2 \abs{1+i(u_j-v_i)}^{-1+\epsilon} du \ll \abs{1+i(v_2-v_1)}^{-1+\epsilon} \abs{1+i(v_3-v_2)}^{-\frac{3}{2}+\epsilon}, \]
	\[ \int_{\R} \abs{1+i(u-v_1)}^{-\frac{3}{4}+\epsilon} \abs{1+i(u-v_2)}^{-\frac{3}{4}+\epsilon} du \ll \abs{1+i(v_2-v_1)}^{-\frac{1}{2}+\epsilon}, \]
	\[ \int_{\R} \prod_{i=1}^3 \abs{1+i(u-v_i)}^{-\frac{3}{2}+\epsilon} du \ll \abs{1+i(v_2-v_1)}^{-\frac{3}{2}+\epsilon} \abs{1+i(v_3-v_2)}^{-\frac{3}{2}+\epsilon}. \]
\end{lem}
The proof matches that of \cite[Lemmas 4 and 6]{Me01}, and we omit it.

The function $E_F(s,t)$ on $t_2 \ge t_1$ satisfies
\[ E_F(s,t) \le E_F(s,t') \]
whenever $t_2' \ge t_2$ and $t_1'+t_2'\ge t_1+t_2$.
With this and the bounds \eqref{eq:Sw4Bound}-\eqref{eq:SwlBound} in mind, bounding the Kloosterman sum side of the Kuznetsov formula is a simple matter of contour shifting.
The polynomial parts of the integrands in the Mellin-Barnes integrals of the various $K^{\pm\pm}_{w_l}$ match, so we consider only the $\varepsilon=(+1,+1)$ term:
Set $\wtilde{F}(\mu) = F(\mu) \cosmu(\mu) \specmu^1(\mu)$, then for $y_1,y_2 > 0$, and $\eta>0$ sufficiently small,
\begin{align*}
	& 4 \abs{y_1 y_2} H_{w_l}(F;y) =\\
	& \pi^{\frac{7}{2}} \sum_{w\in W} \int_{\Re(\mu)=u_4(w)} \wtilde{F}(\mu) p_{\rho+\mu^w}(y) \Lambda_{(-1,-1,-1)}(2 \mu^{w w_l}) d\mu \\
	&+\sum_{w\in W_3} \int_{\Re(\mu)=u_3(w)} \wtilde{F}(\mu) (\pi y_2)^{1+\mu_1^w} \\
	& \qquad \times \int_{\Re(s_1)=-\frac{1}{2}-3\eta} (\pi y_1)^{1-s_1} \Gamma(\mu_2^w-\mu_1^w) \Gamma(\mu_3^w-\mu_1^w) \Gamma(s_1-\mu_2^w) \Gamma(s_1-\mu_3^w) \frac{ds_1}{2\pi i} d\mu \\
	&+\sum_{w\in W_3} \int_{\Re(\mu)=u_2(w)} \wtilde{F}(\mu) (\pi y_1)^{1-\mu_3^w} \\
	& \qquad \times \int_{\Re(s_2)=-\frac{1}{2}-3\eta} (\pi y_2)^{1-s_2} \Gamma(\mu_3^w-\mu_1^w) \Gamma(\mu_3^w-\mu_2^w) \Gamma(s_2+\mu_1) \Gamma(s_2+\mu_2) \frac{ds_2}{2\pi i} d\mu \\
	&+\int_{\Re(\mu)=u_1} \wtilde{F}(\mu) \int_{\Re(s)=(-\frac{1}{2}-\eta,-\frac{1}{2}-\eta)} (\pi y_1)^{1-s_1} (\pi y_2)^{1-s_2} G\paren{0,2s,2\mu} \frac{ds}{(2\pi i)^2} d\mu,
\end{align*}
where we choose
\begin{align*}
	u_1=&(0,0,0), & (u_2(w))^w =& (\tfrac{1}{4}+2\eta,\tfrac{1}{4}+2\eta,-\tfrac{1}{2}-4\eta), \\
	u_3(w) =& -u_2(w w_l), & (u_4(w))^w =& (\tfrac{1}{2}+\eta,0,-\tfrac{1}{2}-\eta).
\end{align*}
Applying Stirling's formula and Lemma \ref{lem:ElemBounds}, we see
\[ \frac{H_{w_l}(F;y)}{(y_1 y_2)^{\frac{1}{2}+\eta}} \ll E_F(u_1,(0,\tfrac{1}{2}))+E_F(u_2(I),(-\tfrac{1}{2},\tfrac{1}{2}))+E_F(u_3(I),(-\tfrac{1}{2},\tfrac{1}{2}))+E_F(u_4(I),(-\tfrac{1}{2},0)), \]
as $y \to 0$, which is sufficient for absolute convergence of the sum of Kloosterman sums.

In a similar manner, we see
\[ H_{w_4}(F;y) \ll y_1^{1-\eta} E_F(0,(-\tfrac{1}{2},\tfrac{1}{2}))+y_1^{\frac{1}{2}+\eta} E_F(u_2(I),(-\tfrac{1}{4},\tfrac{3}{4})), \]
as $y_1 \to 0$.
Of course, we only need $y_1^\epsilon$ for absolute convergence of the sum of Kloosterman sums here, but this gives a better, more concise bound.
This bound could easily be improved by a careful choice of contours in $\mu$ (not vertical lines).

Taking a Hecke eigenbasis for $\mathcal{S}_2^1$ and using the explicit form \eqref{eq:MaxParaRS} below, the maximal parabolic Eisenstein series term is bounded by $E_F^*$, which follows from known lower bounds on the $L$-functions in the denominator \cite{HoffLock,HoffRam}.
One could likely remove the $E_F^*$ term by applying the $GL(2)$ Kuznetsov formula, but this is sufficient for our purposes.

\subsection{The proof of Corollary \ref{cor:Weyl}}
The proof proceeds in two steps.
First, we obtain an upper bound of the correct order of magnitude:
Define
\[ \omega(\mu) = \frac{1-(\mu_1-\mu_2)^2}{100-(\mu_1-\mu_2)^2}, \]
and note that for $\mu = \mu_\varphi$ the Langlands parameters of $\varphi \in \mathcal{S}^1_3$, we have $\omega(\mu) \asymp 1$ by Theorem \ref{thm:KS}.
Taking a test function
\[ F(\mu) = \omega(\mu) \paren{\sum_{w\in W_2} \exp ((\mu_1^w-T\mu_1')^2+(\mu_2^w-T\mu_2')^2)}^2, \]
it follows from Theorem \ref{thm:TechnicalWeyl} that
\begin{align}
\label{eq:RoughUpper}
	\sum_{\norm{\mu_\varphi-T \mu'} < 100} \frac{\abs{\rho_\varphi^*(1)}^2}{\cosmu^1(\mu_\varphi)} \asymp \sum_{\norm{\mu_\varphi-T \mu'} < 100} \frac{\abs{\rho_\varphi^*(1)}^2}{\cosmu^1(\mu_\varphi)} \omega(\mu_\varphi) \ll T^{3+\epsilon}
\end{align}
for $\mu'\in i\mathfrak{a}_\R^*$.
Note that $F(\mu)$ is positive on the spectrum since we have forced the spectral parameters to be of the form
\[ \Re(\mu)=0, \qquad \text{ or } \qquad (x+it,-x+it,-2it) \]
by our choice of Whittaker functions (see the discussion at the end of section 9.3 in \cite{HWII}).
Similarly, the weighted count of complementary series forms $\varphi$ with spectral parameters $\mu_\varphi=(x+it,-x+it,-2it)$, $0\ne \abs{x} \le \frac{1}{2}$, $\abs{t-T} < 100$ is $O(T^{2+\epsilon})$ using a test function of the form
\[ F(\mu) = \omega(\mu) \exp\paren{(\mu_1-\mu_2)^2 + (\mu_1-iT)^2 + (\mu_2-iT)^2}. \]
Note: This is where the need for Theorem \ref{thm:KS} arises, as otherwise, having only the trivial bound $\abs{\Re(\mu_i)} < \frac{1}{2}$, we could not rule out the possibility of a large number of complementary series forms with $\omega(\mu)$ correspondingly small, and so it would be impossible to remove the analytic weights $\omega(\mu)$.

Now, we proceed to the main proof:
Let $\chi_{T\Omega}$ be the characteristic function of the set $T\Omega$, then we define our test function by convolution with a Gaussian approximation to the identity of width $1/\sqrt{\log T}$:
\[ H_T(\mu) = -\frac{\log T}{\pi} \int_{\Re(\mu')=0} \chi_{T\Omega}(\mu-\mu') T^{((\mu_1')^2+(\mu_2')^2)} d\mu', \]
\[ F(\mu) = \omega(\mu) H_T(\mu). \]
Substituting $\mu'\mapsto\mu-\mu'$, this extends to an entire function of $\mu$; by the symmetry of $\Omega$, we have $F(\mu^{w_2})=F(\mu)$.

Let $B(r,\mu)\subset i\mathfrak{a}_\R^*$ be the ball of radius $r$ centered at $\mu$.
At positive distance from the boundary of $T\Omega$, say for $\mu\in i\mathfrak{a}_\R^*$ not in $\mathcal{O}=\partial T\Omega+B(10,0)$, the difference $\omega(\mu) \chi_{T\Omega}(\mu)-F(\mu)$ is negligibly small by the exponential decay of $T^{((\mu_1')^2+(\mu_2')^2)}$.
On the other hand, the assumption on the Minkowski dimension tells us that $\partial \Omega$ can be covered by $O(T^{1+\epsilon})$ balls of radius $\frac{1}{T}$, but simple scaling implies that $\mathcal{O}$ can be covered by $O(T^{1+\epsilon})$ balls of radius $11$, hence
\begin{align}
\label{eq:SmoothWeyl}
	\sum_{\mu_\varphi \in T\Omega} \frac{\abs{\rho_\varphi^*(1)}^2}{\cosmu^1(\mu_\varphi)} =& \sum_{\varphi} F(\mu_\varphi)\frac{\abs{\rho_\varphi^*(1)}^2}{\cosmu^1(\mu_\varphi)}+O\paren{T^{4+\epsilon}}, \\
	\int_{\Re(\mu)=0} F(\mu) \specmu^1(\mu) d\mu =& \int_{T\Omega} \specmu^1(\mu) d\mu+O\paren{T^{4+\epsilon}}, \nonumber
\end{align}
using
\[ \omega(\mu_\varphi) = 1 - \frac{99}{100-(\mu_{\varphi,1}-\mu_{\varphi,2})^2}, \]
and the upper bound \eqref{eq:RoughUpper}.
In order to justify \eqref{eq:SmoothWeyl} on the complementary series spectrum (where $F$ does not approximate $\omega \chi_{T\Omega}$), we needed the bound
\begin{align}
\label{eq:TestFunOffSymmBound}
	F(\mu) \ll& T^{\Re(\mu_1)^2+\Re(\mu_2)^2+\epsilon} \chi_{T\Omega+B(10,0)}(i\Im(\mu))+(\norm{\mu}+T)^{-97}, \quad \text{ on } \quad \abs{\Re(\mu_i)}<1,
\end{align}
which again follows from the rapid decay of $T^{((\mu_1-\mu_1')^2+(\mu_2-\mu_2')^2)}$.
Therefore, the contribution of the complementary series to \eqref{eq:SmoothWeyl} is $O(T^{7/2+\epsilon})$ even using the trivial bound $\abs{\Re(\mu_i)} < \frac{1}{2}$.

Part 1 of the corollary now follows from Theorem \ref{thm:TechnicalWeyl} and \eqref{eq:TestFunOffSymmBound}.
Part 2 is similar.

\section{Rankin-Selberg}
\label{sect:RS}
We recall some facts from \cite[section 6]{HWI} about the spherical, maximal parabolic Eisenstein series induced from the constant function on $GL(2)$.
This function is constructed by \cite[(5.2)]{HWI},
\[ E^0(g,1,s) := \sum_{P_{21}(\Z)\backslash\Gamma} I^0(\gamma g,1,s), \qquad I^0(xyk,1,s) := (y_1^2 y_2)^{\frac{1}{2}+s}. \]
The completion
\[ E^{0*}(g,1,s) := \Gamma_\R(\tfrac{3}{2}+3s) \zeta(\tfrac{3}{2}+3s) E^0(g,1,s), \]
is entire in $s$ except for simple poles at $s=\pm\frac{1}{2}$ with residues $\pm\frac{2}{3}$, and satisfies the functional equation
\[ E^{0*}(g,1,s) = E^{0*}(\vpmpm{--} g^\iota w_l,1,-s). \]
To see this from \cite[Proposition 5.1]{HWI}, note that the Hecke eigenvalues of the constant function $\Phi=1$ on $GL(2)$ are $\lambda_\Phi(n) = \sum_{ab=n} a^{1/2} b^{-1/2}$ so the Hecke $L$-function is
\[ L(\Phi,s)=\zeta(s+\tfrac{1}{2})\zeta(s-\tfrac{1}{2}), \]
or one may find this in \cite[Corollary 2.5]{Fried01}.

\subsection{$L$-functions}
Let us return to our discussion of Theorem \ref{thm:StadesFormula} and Corollary \ref{cor:RSLfunc}; we generally follow \cite{Fried01} here.
If we write the Fourier-Whittaker expansion of Hecke-normalized cusp forms $\varphi$ as
\begin{align}
\label{eq:HeckeFourierWhitt}
	\varphi(g) =& \sum_{m_1=1}^\infty \sum_{m_2 \ne 0} \frac{\lambda_\varphi(m)}{\abs{m_1 m_2}} \sum_{\gamma\in U(\Z)\backslash SL(2,\Z)} W^{1*}(\wtilde{m}\gamma g, \mu_\varphi), \qquad \wtilde{m}=\diag(m_1 m_2,m_1,1)
\end{align}
then the completed $L$-function attached to $\varphi$ (see \cite[page 101]{HIM01} or \cite[appendix A]{MillerSchmid}) is $\Lambda(s,\varphi)=L_\infty(s,\varphi) L(s,\varphi)$ where
\begin{align}
\label{eq:Lvarphi}
	L_\infty(s,\varphi) =& \Gamma_\R(1+s+\mu_{\varphi,1}) \Gamma_\R(1+s+\mu_{\varphi,2}) \Gamma_\R(s+\mu_{\varphi,3}), \quad L(s,\varphi) =& \sum_{n=1}^\infty \frac{\lambda_\varphi(1,n)}{n^s}.
\end{align}

The finite part of the Rankin-Selberg $L$-function attached to Hecke-Maass cusp forms $\varphi$ and $\wbar{\varphi'}$ is
\[ L(s,\varphi \times \wbar{\varphi'}) = \zeta(3s) \sum_{m\in\N^2} \frac{\lambda_\varphi(m) \wbar{\lambda_{\varphi'}(m)}}{(m_1^2 m_2)^s}, \]
and by unfolding (see \cite[Theorem 3.2]{Fried01} and note the extra factor $\frac{1}{2}$ in that definition of the Eisenstein series), we have
\begin{align}
\label{eq:RSIntegral}
	\int_{\Gamma\backslash G} E^{0*}(g,1,s) \varphi(g) \wbar{\trans{\varphi'(g)}} dg = \Lambda(\tfrac{1}{2}+s,\varphi \times \wbar{\varphi'}).
\end{align}
Note that, for $\Re(s)$ in a compact set, $E^{0*}(xyk,1,s)$ is polynomially bounded in the coordinates of $y$, while the cusp forms have exponential decay as $y_i \to \infty$, so the integral defines an entire function in $s$, except for the simple poles of the Eisenstein series.

For a Hecke cusp form $\varphi \in \mathcal{S}^1_3$, it is well-known that $\wbar{\lambda_\varphi(m_1,m_2)}=\lambda_\varphi(m_2,m_1)$, and from \eqref{eq:DualWhitt} it is not too hard to see the Fourier coefficients of the dual form $\wcheck{\varphi}(g) := \varphi(\vpmpm{--}g^\iota w_l)$ are also
\[ \int_{U(\Z)\backslash U(\R)} \wcheck{\varphi}(xg) \wbar{\psi_m(x)} dx=-\frac{\lambda_\varphi(m_2,m_1)}{\abs{m_1 m_2}} W^{1*}(\wtilde{m}g, -\mu). \]
Then substituting $g \mapsto \vpmpm{--} g^\iota w_l$ in \eqref{eq:RSIntegral} gives
\begin{align}
\label{eq:RSIntegral2}
	\Lambda(\tfrac{1}{2}+s,\varphi \times \wbar{\varphi'}) =& \int_{\Gamma\backslash G} E^{0*}(g,1,-s) \wcheck{\varphi}(g) \wbar{\trans{\wcheck{\varphi'}(g)}} dg = \Lambda(\tfrac{1}{2}-s,\wbar{\varphi} \times \varphi'),
\end{align}
which is the functional equation.

The Rankin-Selberg convolution of a Hecke-Maass cusp form with itself factorizes into the exterior- and symmetric-square $L$-functions
\[ L(s,\varphi \times \varphi) = L(s, \ExtSq \varphi) L(s, \SymSq \varphi), \]
and in the case of $GL(3)$, we know (see \cite{Kont01}) the exterior square is actually the dual $L(s, \ExtSq \varphi) = L(s, \wbar{\varphi})$.
Comparing the completions, we have
\begin{align}
\label{eq:LinfSymSq}
	L_\infty(s,\SymSq \varphi) =& \frac{2 \Gamma_\R(3s) \Psi^1(\mu_\varphi,\mu_\varphi,s)}{L_\infty(s,\wbar{\varphi})} \\
	=& \paren{\prod_{i=1}^3 \Gamma_\R(s+2\mu_i)} \Gamma_\R(1+s-\mu_1) \Gamma_\R(1+s-\mu_2) \Gamma_\R(s-\mu_3). \nonumber
\end{align}

The proof of Theorem \ref{thm:KS} given in \cite{KS01} applies functoriality for the holomorphy and functional equations of the symmetric-square $L$-function, which does not assume the cusp form is spherical, and then applies the method of \cite{LRS} with the Rankin-Selberg $L$-function replaced with the symmetric-square $L$-function.
Although the theorem given in \cite{LRS} only applies to spherical forms, the formula \eqref{eq:LinfSymSq} shows that the argument still applies to the weight-one forms as the gamma factors have some of the poles shifted farther to the \textit{left}, but the relevant factors
\[ \Gamma_\R\paren{s+2\mu_1} \Gamma_\R\paren{s+2\mu_2} \Gamma_\R\paren{s+2\mu_3} \]
remain the same, so the argument goes through verbatim.

In fact, the results of \cite{LRS} extend to non-spherical forms on $GL(n)$ (and \cite{KS01} up to $GL(4)$) as follows:
The appendix to \cite{MillerSchmid} very nicely writes out the relevant gamma factors; in that notation, for a cusp form $\varphi$ coming from an induced representation
\[ I(P; \sigma_1[s_1],\ldots,\sigma_r[s_r]) \]
on $GL(n,\R)$ (over $\Q$) with each $\sigma_i$ one of $D_k$, $\sgn$, or $\mathrm{triv}$, the gamma factors of the Rankin-Selberg convolution $L_\infty(s', \varphi \times \varphi)$ contain the factors
\[ L_\infty(s', \sigma_i[s_i] \times \sigma_i[s_i]) = \piecewise{\Gamma_\R(s'+2s_i) & \If \sigma_i \in\set{\mathrm{triv}, \sgn}, \\ \Gamma_\C(s'+2s_i) & \If \sigma_i = D_k,} \]
and each of these contains a gamma factor $\Gamma_\R(s'+2s_i)$.
Note that unitaricity implies that $\sigma_i[-\wbar{s_i}]$ is also a component.
So we conclude as in \cite{LRS} that $\abs{\Re(s_i)} \le \frac{1}{2}-\frac{1}{n^2+1}$ (and $\abs{\Re(s_i)} \le \frac{9}{22}$ on $GL(4)$ by \cite{KS01}).

\subsection{The normalization coefficients}
For a Hecke-normalized Maass cusp form $\varphi$ as in \eqref{eq:HeckeFourierWhitt}, we have the factorization
\begin{align}
\label{eq:AdSqFactor}
	L(s,\varphi \times \wbar{\varphi}) = \zeta(s) L(s,\AdSq \varphi),
\end{align}
and taking the residue of \eqref{eq:RSIntegral} at $s=\frac{1}{2}$ when $\varphi'=\varphi$ gives
\[ \frac{2}{3} \norm{\varphi}^2 = \frac{L(1,\AdSq \varphi)}{\pi \cosmu^1(\mu_\varphi)}. \]
This implies
\begin{align}
	\frac{\wbar{\rho_\varphi^*(m)} \rho_\varphi^*(n)}{\cosmu^1(\mu_\varphi)} =& \frac{2\pi}{3} \frac{\wbar{\lambda_\varphi(m)} \lambda_\varphi(n)}{L(1,\AdSq \varphi)}
\end{align}
for an $L^2$-normalized Hecke-Maass cusp form $\varphi$.

From \cite[(3.31),(5.6),(5.19)]{HWI}, \cite[Theorem 5]{HWII} and
\[ \Lambda_{(0,1,1)}(\mu) \Lambda_{(0,1,1)}(-\mu) = \frac{1}{\pi \cosmu^1(\mu)}, \]
we have
\begin{align}
	\frac{\wbar{\rho_\phi^*(m;\mu_1)} \rho_\phi^*(n;\mu_1)}{\cosmu^1(\mu_1+\mu_\phi,\mu_1-\mu_\phi,-2\mu_1)} =& 4\pi \frac{\wbar{\lambda_\phi(m,\mu_1)} \lambda_\phi(n,\mu_1)}{L(\phi,1+3\mu_1) L(\phi,1-3\mu_1)},
\end{align}
when $\phi$ is Hecke-normalized.
When $\phi$ is $L^2$-normalized (see \cite[section 5.3]{HWI}), this becomes
\begin{align}
\label{eq:MaxParaRS}
	\frac{\wbar{\rho_\phi^*(m;\mu_1)} \rho_\phi^*(n;\mu_1)}{\cosmu^1(\mu_1+\mu_\phi,\mu_1-\mu_\phi,-2\mu_1)} =& 2\pi \frac{\wbar{\lambda_\phi(m,\mu_1)} \lambda_\phi(n,\mu_1)}{L(\phi,1+3\mu_1) L(\phi,1-3\mu_1) L(1,\AdSq \phi)}.
\end{align}

\bibliographystyle{amsplain}

\bibliography{HigherWeight}

\end{document}